\documentclass[11pt]{article}
\usepackage{amssymb,amsmath,amsthm}
\usepackage{graphicx}

\newtheorem{thm}{Theorem}
\newtheorem*{theo}{Theorem}
\newtheorem*{propo}{Proposition}

\newtheorem{lemma}{Lemma}

\newtheorem{lemm}{Lemma}[section]
\newtheorem{corol}[lemm]{Corollary}
\newtheorem{defi}[lemm]{Definition}
\newtheorem{quest}[lemm]{Question}
\newtheorem{exer}[lemm]{Exercise}
\newtheorem{exam}[lemm]{Example}
\newtheorem{rema}[lemm]{Remark}
\newtheorem{remas}[lemm]{Remarks}
\newtheorem{propos}[lemm]{Proposition}
\newtheorem{notas}[lemm]{Notations}
\newtheorem{notat}[lemm]{Notation}

\def\demo{\medskip\goodbreak\noindent
     \hbox{\sc Proof \kern .3em}\ignorespaces}%
     \def\idemo{\medskip\goodbreak\noindent
     \hbox{\sc Ideas of proof \kern .3em}\ignorespaces}%
  \def \qedbox{$\square$}%
  \def \qed{\hglue1mm\hfill{\ifmmode\qedbox
     \else\unskip\ \hglue0mm\hfill\qedbox\medskip
      \goodbreak\fi}}%
\def\enddemo{\qed\goodbreak\vskip10pt}%

\def\qed{\hglue1mm\hfill\raise -2pt\hbox{\vrule\vbox to 10pt{\hrule width
4pt
                  \vfill\hrule}\vrule}}

\newcommand{\T}{\mathbb {T}}

\newcommand{\A}{\mathbb {A}}

\newcommand{\R}{\mathbb {R}}
\newcommand{\Q}{\mathbb {Q}}

\newcommand{\Z}{\mathbb {Z}}
\newcommand{\N}{\mathbb {N}}

\newcommand{\Vc}{\mathcal {V}}
\newcommand{\Uc}{\mathcal {U}}

\newcommand{\Dc}{\mathcal {D}}
\newcommand{\Wc}{\mathcal {W}}

\newcommand{\Cc}{\mathcal {C}}
\newcommand{\Ic}{\mathcal {I}}

\newcommand{\Mc}{\mathcal {M}}
\newcommand{\Kc}{\mathcal {K}}
\newcommand{\Ec}{\mathcal {E}}
\newcommand{\Gc}{\mathcal {G}}

\newcommand{\Fc}{\mathcal {F}}
\newcommand{\Ac}{\mathcal {A}}
\newcommand{\Sc}{\mathcal {S}}

\begin{document}
\title{Hyperbolicity for conservative twist maps of the 2-dimensional annulus}
\author{M.-C. ARNAUD \date{}
\thanks{ANR-12-BLAN-WKBHJ}
\thanks{Avignon Universit\'e , Laboratoire de Math\'ematiques d'Avignon (EA 2151),  F-84 018 Avignon,
France. e-mail: Marie-Claude.Arnaud@univ-avignon.fr} 
\thanks{membre de l'Institut universitaire de France}}
\maketitle
\abstract{These are notes for a minicourse given at Regional Norte UdelaR in Salto, Uruguay for
the conference ÒCIMPA Research School - Hamiltonian and Lagrangian DynamicsÓ.
We will present Birkhoff and Aubry-Mather theory for the conservative twist maps of the 2-dimensional annulus and  focus on what happens close to the Aubry-Mather sets: definition of the Green bundles, link between hyperbolicity and shape of the Aubry-Mather sets, behaviour close to the boundaries of the instability zones.  We will also give some open questions.  This course is the  second part of a minicourse that was begun by R.~Potrie. Some topics of the part of R.~Potrie will be useful for this part.

Many thanks to E. Maderna and L. Rifford for the invitation to give the mini-course and to R.~Potrie for
accepting to share the course with me.}

 \vfill\eject

\tableofcontents

  \newpage
\section{Introduction to conservative twist maps}
\begin{notas}{\rm
\begin{enumerate}
\item[$\bullet$] $\T=\R/\Z$ is the circle; $\A=\T\times \R$ is the annulus and $(\theta, r)\in\A$ refers to a point of $\A$;
\item[$\bullet$] $\A$ is endowed with its symplectic form $\omega=dr \wedge d\theta= d\lambda$ where $\lambda=rd\theta$ is the Liouville 1-form;
\item[$\bullet$] $p:\R^2\rightarrow \A$ is the universal covering;
\item[$\bullet$] $\pi:\A\rightarrow \T$ is the first projection: $\pi(\theta , r)=\theta$ and $\pi:\R^2\rightarrow \R$ is its lift, which is also a projection: $\pi(\theta , r)=\theta$;
\item[$\bullet$]  for every point $x=(\theta,r)$, the  vertical  line at $x$  is $\Vc(x)=\{ \theta\}\times \R\subset \R^2$ or $\Vc(x)=\{ \theta\}\times \R\subset \A$;
\item[$\bullet$] the vertical subspace is the tangent subspace to the vertical line: $V(x)=T_x \Vc(x)$;
\item[$\bullet$]  all the measures we will deal with are assumed to be Borel probabilities. The support of $\mu$ is denoted by ${\rm supp}\mu$.
\end{enumerate}}
\end{notas}
 If $x\in M$ is an elliptic periodic point of a Hamiltonian flow that is defined on a 4-dimensional symplectic manifold $M$, using symplectic polar coordinates in an annular Poincar\'e section contained  in the  energy level of $x$, we obtain in general a  first return map $T:\Ac\rightarrow \A$  that is defined on some bounded sub-annulus $\Ac$  of $\A$ by $T(\theta, r)=(\theta +\alpha+\beta r, r)+o(r)$ with $\beta\not=0$. This is locally a {\em conservative twist map}. 
  \begin{center}
\includegraphics[width=8cm]{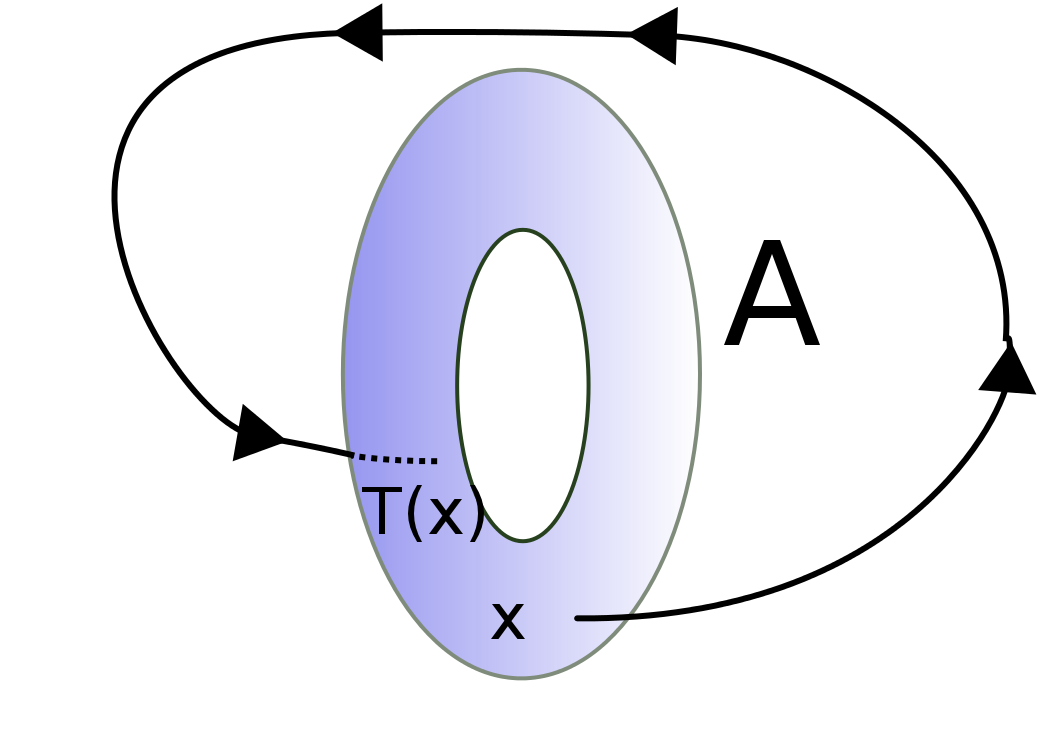}\end{center} 
\begin{defi} {\rm A {\em positive (resp. negative) twist map} is a $C^1$-diffeomorphism $f:\A\rightarrow \A$ such that
\begin{enumerate}
\item $f$ is isotopic to the identity map ${\rm Id}_\A$ (i.e. $f$ preserve the orientation and the two boundaries of the annulus);
\item $f$ satisfies the {\em twist condition} i.e. there exists $\varepsilon>0$ such that  for any $x\in \A$, we have: $\frac{1}{\varepsilon}>D(\pi\circ f)(x)(0, 1)>\varepsilon$ (resp. $-\frac{1}{\varepsilon}<D(\pi\circ f)(x)(0, 1)<-\varepsilon$). In the first case the twist is positive, in the second case it is negative.
\begin{center}
\includegraphics[width=6cm]{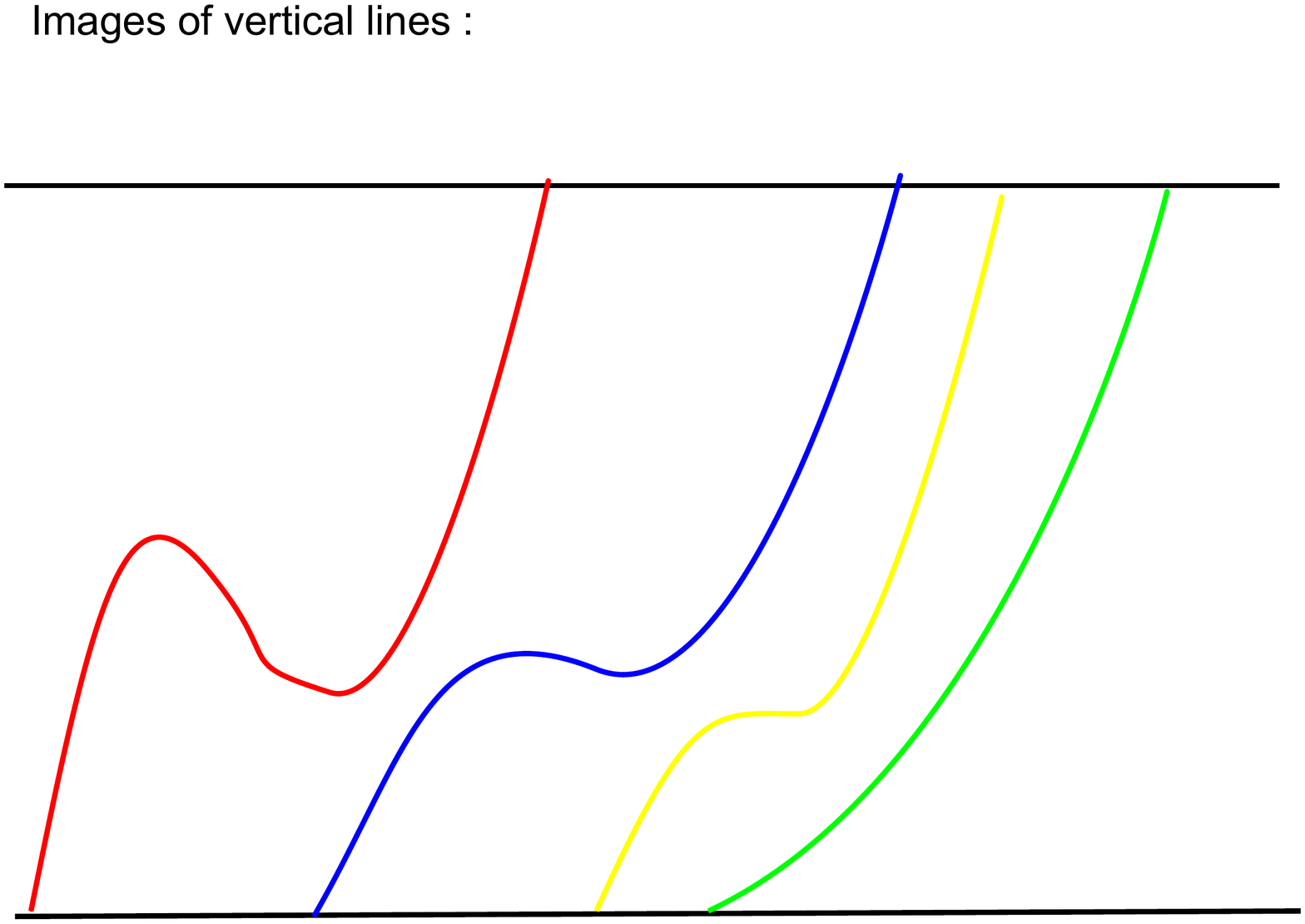}\end{center} 
The twist map is {\em conservative } (or {\em exact symplectic}) is $f^*\lambda-\lambda$ is an exact 1-form.
\end{enumerate}}

\end{defi}

\begin{remas}{\rm  \begin{enumerate}
\item Saying that the diffeomorphism $f$ is isotopic to identity means that:
\begin{enumerate}
\item[$\bullet$] $f$ preserves the orientation;
\item[$\bullet$] $f$ fixes the two ends $\T\times \{-\infty\}$ and $\T\times\{+\infty\}$ of the annulus.
\end{enumerate}
\item The reader can ask why we don't just ask that $f$ preserves the area form (symplectic form) $\omega$, i.e. $0=f^*\omega-\omega=d(f^*\lambda-\lambda)$. We ask not only that $f^*\lambda-\lambda$ is closed, we ask that it is exact. Indeed, we want to avoid symplectic twist maps as $(\theta, r)\mapsto (\theta+r, r+1)$: all the orbits come from $\T\times \{ -\infty\}$ and go to $\T\times \{ +\infty\}$ and there is no non-empty compact invariant set for such a map. We will see in section \ref{SAubryMather} that this never happens for exact symplectic twist maps;
\item Note that $f$ is a positive conservative twist map if and only if $f^{-1}$ is a negative conservative twist map. Hence from now we will assume that all the considered conservative twist maps are positive.
\end{enumerate}}
\end{remas}

\begin{exer}
{\rm Let $f:\A\rightarrow\A$ be a conservative twist map. Using Stokes formula, prove that if $\gamma:\T\rightarrow \A$ is a $C^1$-embedding, then the (algebraic) area of the domain that is  between $\gamma$ and $f(\gamma)$ is zero.\begin{center}
\includegraphics[width=4cm]{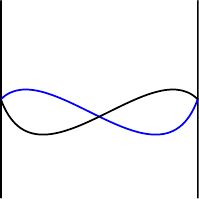}\end{center} }

\end{exer}

\begin{exam}\label{Ex1}
{\rm Consider the map we introduced by using polar coordinates for a first return map $T(\theta, r)=(\theta+\alpha+\beta r,r)$ and   assume that $\beta>0$ (or replace $T$ by $T^{-1}$). Then $D(\pi\circ T)\begin{pmatrix} 0\\ 1\end{pmatrix}=\beta>0$ hence $T$ is a (positive) twist map.  Moreover, $T^*(rd\theta)-rd\theta=\beta rdr=d\left(\frac{\beta}{2}r^2\right)$ hence $T$ is a conservative twist map.

  \begin{center}
\includegraphics[width=4cm]{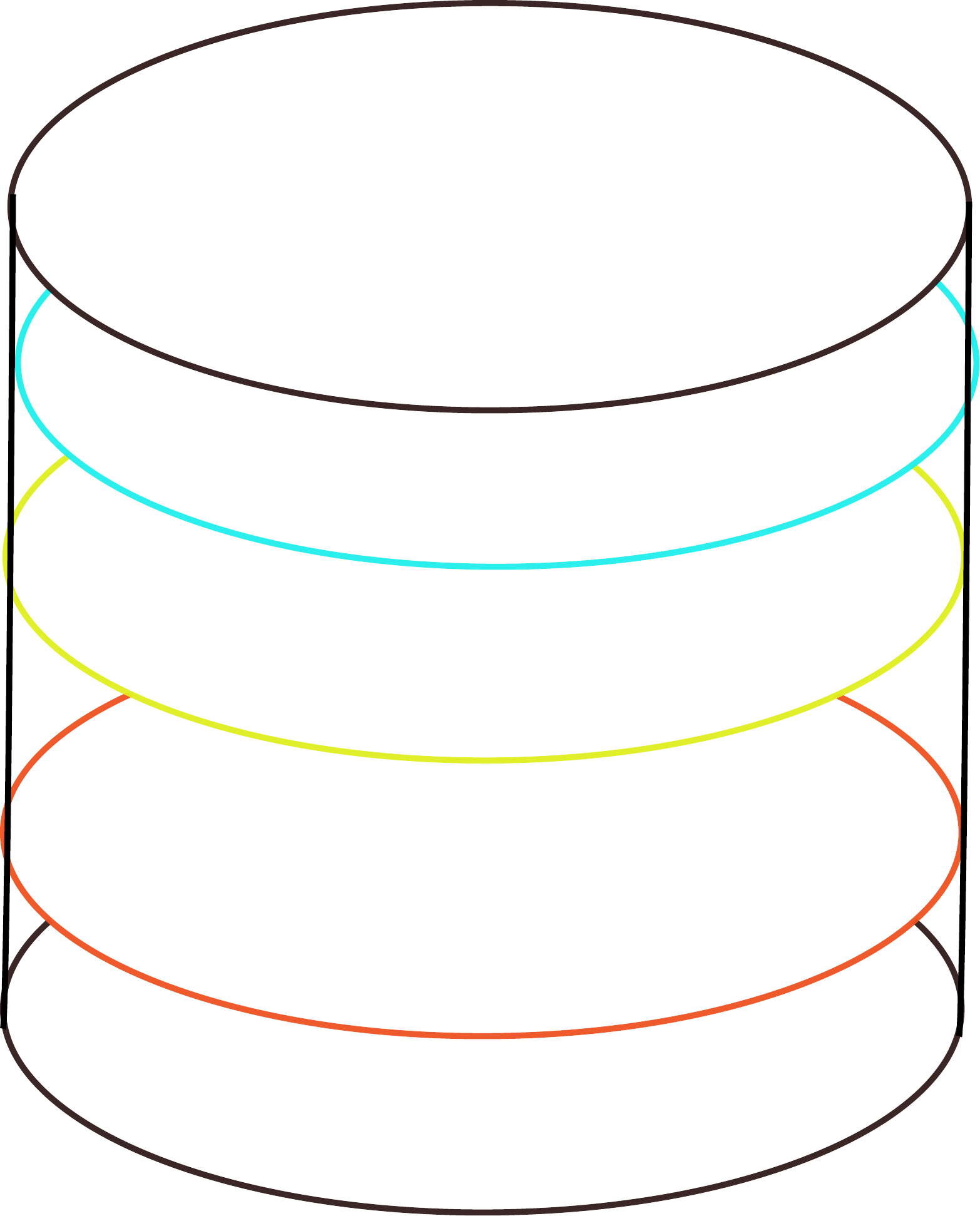}\end{center} 
Note that the dynamics is very simple: the annulus is foliated by invariant circles $\T\times\{ r\}$ and the restriction of $T$ to every such circle is a rotation.}
\end{exam}

\begin{exam}{\rm
The  {\em standard family } depends on a parameter $\lambda\in\R$. It is defined by 
$$f_\lambda(\theta, r)=(\theta+r+\lambda \sin 2\pi \theta, r+\lambda \sin 2\pi\theta).$$

Note that for $\lambda=0$, the map is just the map $T=f_0$ of Example \ref{Ex1}. When $\lambda$ increases from $0$ to $+\infty$, we observe fewer and fewer  invariant graphs.  \begin{center}
\includegraphics[width=5cm]{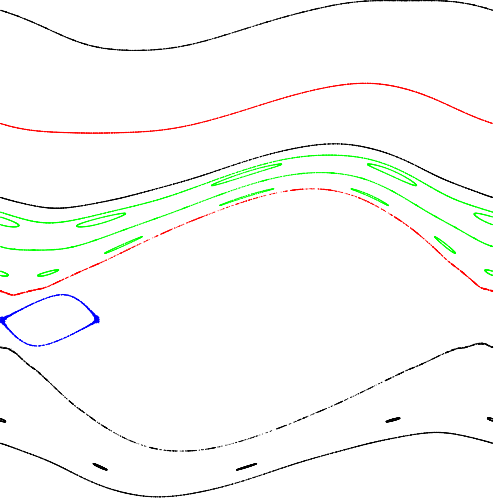}\end{center} 
J.~Mather and S.~Aubry even proved that for  $2\pi\lambda > 4/3$, $f_\lambda$ has no continuous invariant graph. }

\end{exam}
\begin{exer}\begin{enumerate}
\item Check that the functions $f_\lambda$ are all conservative twist maps. \\
Assume that the graph of a continuous map $\psi:\T\rightarrow
 \R$ is invariant by a map $f_\lambda$.
\item Prove that $g_\lambda(\theta) =\theta+\lambda \sin(2\pi\theta)+\psi(\theta)$ is an orientation preserving homeomorphism of $\T$.\\
{\em Hint}: note that $\pi\circ f_\lambda (\theta, \psi(\theta))=g_\lambda(\theta)$.
\item Prove that $g_\lambda^{-1}(\theta)=\theta-\psi(\theta)$.\\
{\em Hint}: prove that $f^{-1}(\theta, r)=(\theta-r, r-\lambda\sin2\pi(\theta-r))$.
\item Check  that $g_\lambda(\theta)+g_\lambda^{-1}(\theta)=2\theta +\lambda\sin2\pi\theta$. Deduce that for $\lambda>\frac{1}{\pi}$, $f_\lambda$ has no continuous invariant graph.
\end{enumerate}
\end{exer}
We can characterize the conservative twist maps by their {\em generating functions}.

\begin{propos}\label{Pgene} Let $F:\R^2\rightarrow \R^2$ be a $C^1$ map. Then $F$ is a lift of a conservative twist map $f:¥A\rightarrow \A$  if and only if there exists a $C^2$ function  such that
\begin{enumerate}
\item[$\bullet$] $\forall \theta, \Theta\in\R, S(\theta+1, \Theta +1)=S(\theta, \Theta)$;
\item[$\bullet$] there exists $\varepsilon>0$ so that for all $\theta, \Theta\in\R$, we have 
$$\varepsilon<-\frac{\partial^2S}{\partial \theta \partial \Theta}(\theta, \Theta)<\frac{1}{\varepsilon};$$
\item[$\bullet$] $F(\theta, r)=(\Theta, R)\Longleftrightarrow R=\frac{\partial S}{\partial \Theta}(\theta, \Theta)\quad{\rm and}\quad r=-\frac{\partial S}{\partial \theta}(\theta, \Theta)$.
\end{enumerate}
\end{propos} 

In this case, we say that $S$ is a {\em generating function} for $F$ (or $f$).
The proof of Proposition \ref{Pgene}  is given in subsection \ref{sstwistgene}.

\begin{exer}{\rm
Check that a generating function of the standard map $f_\lambda$ is $S_\lambda(\theta, \Theta)=\frac{1}{2}(\Theta-\theta)^2-\frac{\lambda}{2\pi}\cos 2\pi\theta$.}
\end{exer}
\begin{rema}{\rm Generating functions are very useful to construct  new examples or perturbations of known examples of conservative twist maps. Indeed, we only need a function to define a 2-dimensional conservative twist map.\\
Using generating functions, we can for example prove that for every $k\in [1, \infty]$, there is a dense $G_\delta$ subset $\Gc$ of the set of $C^k$ conservative twist maps such that at every periodic point $x$ of $f\in \Gc$ with period $n$, $Df^n(x)$ has two distinct eigenvalues (and then these eigenvalues are different from $\pm1$). A similar dense $G_\delta$ subset $\Gc$ exists such that the intersections of the stable and unstable submanifolds of every pair of periodic hyperbolic points transversely intersect (when they intersect).}
\end{rema}

\section{The invariant curves}

\subsection{Invariant continuous graphs and first Birkhoff theorem}

In the '20s,  G.~D.~Birkhoff proved (see \cite{Bir1}) that the invariant continuous graphs by a twist map are locally uniformly Lipschitz.

\begin{thm}\label{Tbir1}{\bf (G.~D.~Birkhoff)}
Let $f:\A \rightarrow \A$ be a conservative twist map and $x\in \A$. Then there exists a $C^1$-neighborhood $\Uc$ of $f$, a neighborhood $U$ of $x$ in $\A$ and a constant $C>0$ such that if the graph of a continuous map $\psi:\T\rightarrow \R$ meets $U$ and is invariant by a $g\in\Uc$, then $\psi$ is $C$-Lipschitz.
\end{thm}
Theorem \ref{Tbir1} is a consequence of a result that concerns all the Aubry-Mather sets and that we will prove later: Proposition \ref{PAMLip}.
\begin{corol} \label{Cbir1}
Let $f:\A \rightarrow \A$ be a conservative twist map and let $K\subset  \A$ be a compact subset of $\A$. Then there exists a $C^1$-neighborhood $\Uc$ of $f$  and a constant $C>0$ such that if the graph of a continuous map $\psi:\T\rightarrow \R$ meets $K$ and is invariant by a $g\in\Uc$, then $\psi$ is $C$-Lipschitz.
\end{corol}

\begin{exer}{\rm 
Prove Corollary \ref{Cbir1}.}
\end{exer}

From Theorem \ref{Tbir1} and Ascoli theorem, we deduce

\begin{corol}\label{Cclosedcurves}
Let $f$ be a conservative twist map of $\A$. The the union $\Ic(f)$ of all its invariant continuous graphs is a closed invariant subset of $f$.
\end{corol}
\begin{exer}{\rm
Prove Corollary \ref{Cclosedcurves}.}
\end{exer}

\begin{remas}{\rm\begin{enumerate}
\item The set $\Ic(f)$ can be empty: this is the case for the standard map $f_{\lambda}$ with $\lambda>\frac{2}{3\pi}$.
\item Using the connecting lemma that was proved by S.~Hayashi in 2006 (see \cite{Hay1}) and more specifically some related results that are contained in \cite{ABC},  Marie Girard proved (in her non-published PhD thesis) that there is dense $G_\delta$ subset $\Gc$ of the set of $C^1$ conservative twist maps such that every $f\in\Gc$ has no continuous invariant graph.
\item Don't deduce that having an invariant graph rarely happens for the conservative twist maps: it depends on their regularity ($C^1, C^3, \dots, C^\infty$). Indeed, the famous theorems K.A.M. (for Kolmogorov-Arnol'd-Moser, see \cite{Arno}, \cite{Kol}, \cite{Rus}) tell us that if a $C^\infty$ conservative twist map $f$ has a $C^\infty$ invariant graph $\Cc$ such that the restriction $f_{|\Cc}$ is $C^\infty$ conjugated to a Diophantine rotation $\theta\mapsto \theta +\alpha$ (i.e. $\alpha$ is Diophantine: there exist $\gamma, \delta>0$ so that for every $p\in\Z$ and $q\in\N^*$, we have $|\alpha-\frac{p}{q}|\geq \frac{\gamma}{q^{1+\delta}}$), there exists a neighborhood $\Uc$ of $f$ in $C^\infty$-topology such that every $g\in\Uc$ has a $C^\infty$ invariant graph $\Gamma$ such that $g_{|\Gamma}$ is $C^\infty$-conjugated to $f_{|\Cc}$.

As the completely integrable standard map $f_0$ has a lot of such invariant graphs, we deduce that for $\lambda$ small enough, $f_\lambda$ has many $C^\infty$ invariant graphs.
\end{enumerate}}\end{remas}

\begin{rema}{\rm We  will see that even when a conservative twist map has no continuous invariant graph, it has a lot of compact invariant subsets: periodic orbits, and even invariant Cantor sets (these are the Aubry-Mather sets, see section \ref{SAubryMather}).}
\end{rema}
%PRLER DYNAMIQUE GENERIQUE DS I(F) (voir MH).

\subsection{Circle homeomorphisms and dynamics on $\Ic(f)$}\label{sscircle}

Now let us explain how is the dynamics restricted to $\Ic(f)$.\\
The dynamics restricted to every invariant graph is Lipschitz conjugated (via $\pi$) to an orientation preserving bi-Lipschitz homeomorphism of $\T$. The classification of the orientation preserving homeomorphisms of the circle is due to H.~Poincar\'e and given in \cite{K-H} (see \cite{He1} for more results). Let us recall quickly the main results. We assume that $h:\T\rightarrow \T$ is an orientation preserving homeomorphism and that $H_1, H_2 :\R\rightarrow \R$ are some lifts of $h$ (then $H_2-H_1=k$ is an integer). Then  
\begin{enumerate}
\item[$\bullet$] the sequence $\left( \frac{H_i^n-Id}{n}\right)_{n\in\N}$ uniformly converge to  a real number  $\rho(H_i)$ that is called the {\em rotation number } of $H_i$; note that $\rho(H_2)-\rho(H_1)=k$; then the class  of $\rho(H_i)$ modulo $\Z$ defines a unique number $\rho(h)\in\T$ and is called the rotation number of $h$;
\item[$\bullet$] $\rho(H_i)=\frac{m}{n}\in\Q$ (with $m$ and $n$ relatively prime) if and only if there exists a point $t\in\R$ so that $H_i^n(t)=t+m$; in this case a point $t$ of $\T$ is either periodic for $h$ or such that there exist two periodic points $t_-$, $t_+$ with period $n$ for $h$ such that $$\displaystyle{\lim_{\ell\rightarrow+\infty} d(h^{-\ell}t, h^{-\ell}t_-)=\lim_{\ell\rightarrow+\infty} d(h^{\ell}t, h^{\ell}t_+)=0}.$$ In this last case, $t$ is negatively heteroclinic to $t_-$ and positively heteroclinic to $t_+$.
\item[$\bullet$] when $\rho(h)\notin \Q/\Z$, $h$ has no periodic points and either the dynamics is minimal and $C^0$-conjugated to the rotation $t\mapsto t+\rho(h)$ or the non wandering set of $h$ is a Cantor subset (i.e.   non-empty compact totally disconnected with no isolated point) $\Omega$, $h_{|\Omega}$ is minimal and all the orbits in $\T\backslash \Omega$ are wandering and homoclinic to $\Omega$ (this means that $\displaystyle{\lim_{\ell\rightarrow \pm\infty}d(h^\ell t, \Omega)=0}$). Moreover, $f$ has a unique invariant measure, and its support is $\Omega$.
\end{enumerate}

 Moreover, if $q\in\Z^*$, $p\in\Z$ are such that $\rho(H_i)<\frac{p}{q}$ (resp. $\rho(H_i)>\frac{p}{q}$), then we have $ H_i^q(t)-t-p<0$ (resp. $ H_i^q(t)-t-p>0$). We deduce that $$\forall k\in\Z, |H_i^k(t)-t-k\rho(H_i)|\leq 1.$$
 
 \begin{defi}
 {\rm When an invariant graph has an {  irrational} (resp. {  rational}) rotation number, we will say that the graph is {\em irrational} (resp. {\em rational}).\\
 When the rotation number is irrational and the dynamics is not minimal, we have a {\em Denjoy counter-example}.
 
 }
 \end{defi}
 
\subsection{Lyapunov exponents of the invariant curves}\label{ssLyapcurves}
\begin{defi}{\rm 
Let $\Cc\subset \A$ be a set that is invariant by a map $f:\A\rightarrow \A$. Then its stable and unstable sets are defined by
$$ W^s (\Cc, f)=\{ x\in \A; \lim_{k\rightarrow +\infty} d(f^kx, \Cc)=0\}$$
and 
$$ W^u (\Cc, f)=\{ x\in \A; \lim_{k\rightarrow +\infty} d(f^{-k}x, \Cc)=0\}.$$ One of these two sets is {\em trivial} if it is equal to $\Cc$.}
\end{defi}
\begin{exam}  {\rm We consider the Hamiltonian flow of the pendulum. In other words, we define $H:\A\rightarrow \R$ by $H(\theta, r)=\frac{1}{2}r^2+\cos2\pi \theta$ and its Hamiltonian flow $(\varphi_t)$ is determined by the Hamilton equations: $\dot\theta=\frac{\partial H}{\partial r}=r$ and $\dot r=-\frac{\partial H}{\partial \theta}=2\pi\sin 2\pi\theta$. For $t>0$ small enough, the time $t$ map $f=\varphi_t$ is a conservative twist map, and as $H$ is constant along the orbits we can find a lot of invariant curves.
\begin{center}
\includegraphics[width=3cm]{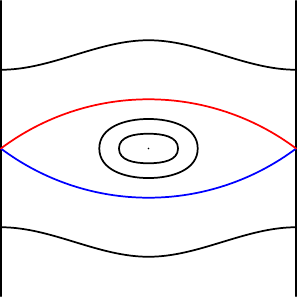}
\end{center}
Note on this picture that there exists two Lipschitz but non $C^1$ invariant graphs, that are the separatrices of the hyperbolic fixed point.\\
Such a  separatrix carries only one invariant ergodic measure, the Dirac mass at the hyperbolic fixed point, and then the Lyapunov exponents of this measure are non zero, and there are non-trivial stable and unstable sets for this separatrix (that is the union of the two separatrices).}
\end{exam}
Hence this is an example of a rational invariant graph that carries an hyperbolic invariant measure. What happens in the irrational case? It is not hard to prove that if the graph of a $C^1$-map is invariant by a conservative twist map and irrational, then the unique ergodic measure supported in the curve has zero Lyapunov exponents. When the invariant curve is just assumed to be Lipschitz, this is less easy to prove but also true as we will see in Theorem \ref{Tirrational}.

\begin{rema} {\rm There exist examples of $C^2$ conservative twist maps that have an irrational invariant Lipschitz graph that is not $C^1$. Such an example is built in \cite{Arna2}. We don't know if such an example exists when the twist map in $C^\infty$ or when the dynamics restricted to the graph is not Denjoy (i.e. has a dense orbit).}
\end{rema}

\begin{quest}{\rm 
Does there exist a $C^\infty$ conservative twist map that has an invariant continuous graph on which the dynamics is Denjoy?}
\end{quest}

\begin{quest}{\rm 
Does there exist a $C^\infty$ conservative twist map that has an invariant irrational continuous graph that is not $C^1$?}
\end{quest}

\begin{quest}{\rm 
If a conservative twist map has an invariant irrational continuous graph on which the restricted dynamics has a dense orbit, is the invariant curve necessarily $C^1$?}
\end{quest}

\begin{remas}{\rm 
\begin{enumerate}
\item From Theorem \ref{Tirrational} and Theorem \ref{TLyapshape} that we will prove later, it is not hard to deduce that if a conservative twist map has an invariant irrational graph $\gamma$  that carries the invariant probability measure $\mu$, then $\gamma$ is $C^1$-regular $ \mu$-almost everywhere (see Definition \ref{DefC1}).
\item In fact, I proved in \cite{Arna1} that any graph that is invariant by a conservative twist map is $C^1$ above a $G_\delta$ subset of $\T$ that has full Lebesgue measure.
\end{enumerate}
}
\end{remas}

With P. Berger, we proved the following result (see \cite{ArnaBer}).

\begin{thm}\label{Tirrational} {\bf (M.-C.~Arnaud \& P.~Berger)}
Let $\gamma$ be an irrational  invariant graph by a $C^{1+\alpha}$ conservative twist map. Then the Lyapunov exponents of the unique invariant probability with support in $\gamma$ are zero. Hence
$$\forall \varepsilon>0, \forall x\in W^s(\gamma, f)\backslash \gamma, \lim_{n\rightarrow +\infty} e^{n\varepsilon} d(f^nx, \gamma)=+\infty.$$
\end{thm}
The convergence to an irrational invariant curve is slower than exponential. We will explain in subsection \ref{ssinstability} that a lot of conservative twist maps have an irrational invariant curve with a non trivial stable set.

\demo
We begin by proving  the first part of the theorem. 

Assume that $\gamma$ is an invariant continuous graph  by a $C^{1+\alpha}$ conservative twist map  $f$ and that some ergodic  invariant probability $\mu$  with support in $\gamma$ is { hyperbolic}, i.e. has two Lyapunov exponents such that $\lambda_1<0<\lambda_2$. As $f$ is symplectic, then $\lambda_2=-\lambda_1=\lambda$.

We use Pesin theory and Lyapunov charts (rectangles $R(f^kx)$) along a generic orbit $(f^kx)$ for $\mu$: in such a chart, the dynamics is almost linear and hyperbolic
\begin{center}
\includegraphics[width=10cm]{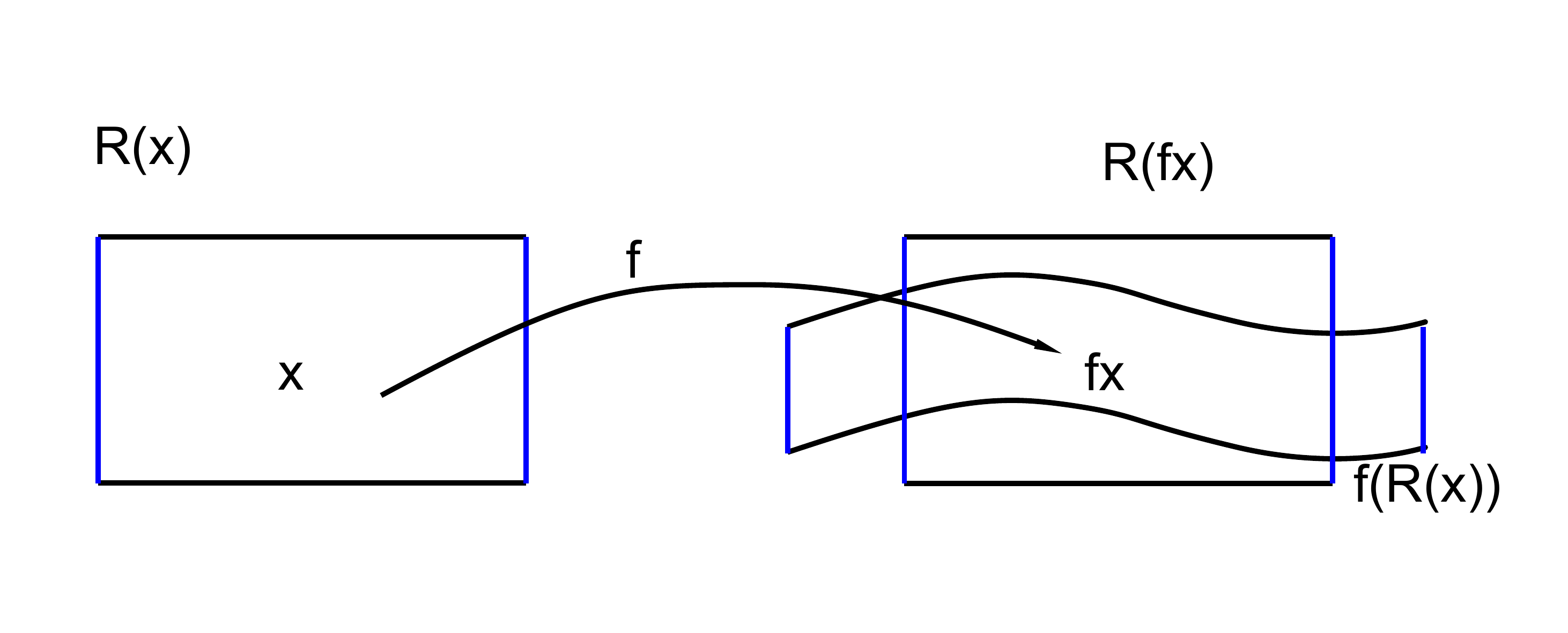}
\end{center} 
We will prove that $\mu$-almost $x$ is periodic. 
The curve $\gamma$ is endowed with some orientation. Note that $f_{|\gamma}$ is orientation preserving.
 
We decompose the boundary $\partial R $ of the domain of a chart   $R $ into $\partial^s R =\{-\rho, \rho\}\times [-\rho, \rho]$ and $\partial^u R =[-\rho, \rho]\times \{ -\rho, \rho\}$
 
\begin{center}
\includegraphics[width=5.6cm]{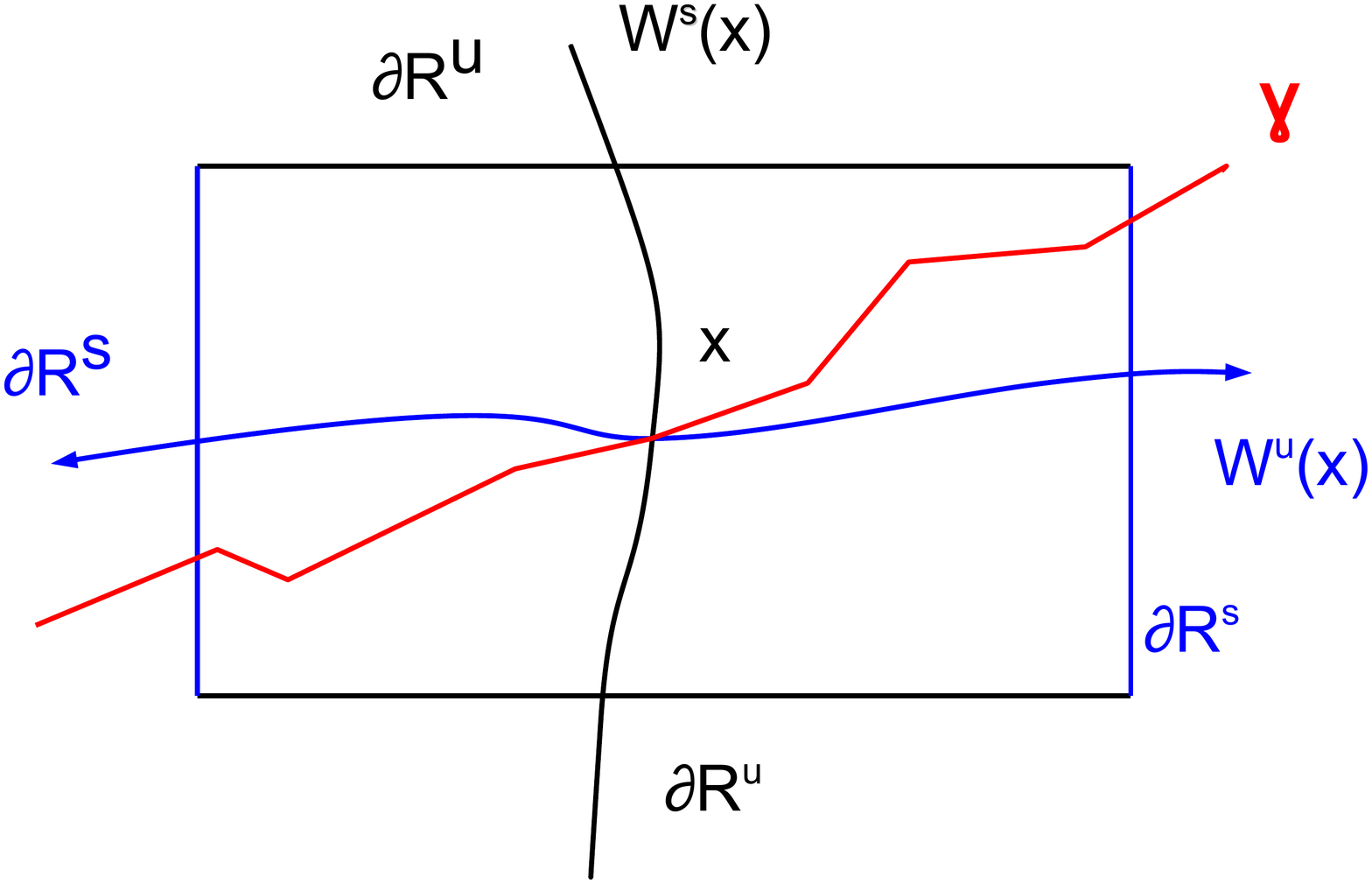}
\end{center} 
 
Let  $\gamma_x$ be the connected components of $\gamma\cap R(x)$ that contains $x$ and let $\eta_x$ be the set of the points of $\gamma_x$ that are after $x$ (for the orientation of $\gamma_x$).

We will prove that $\mu$-almost $x$ is periodic and $\eta_x\subset W^s(x)$ or $\eta_x\subset W^u(x)$.

\begin{lemm}We have either for $\mu$ almost every $x$, $\eta_x(1)\in \partial R^s(x)$ or for $\mu$ almost every $x$, $\eta_x(1)\notin \partial R^s(x)$.  
\end{lemm}
\begin{center}
\includegraphics[width=9cm]{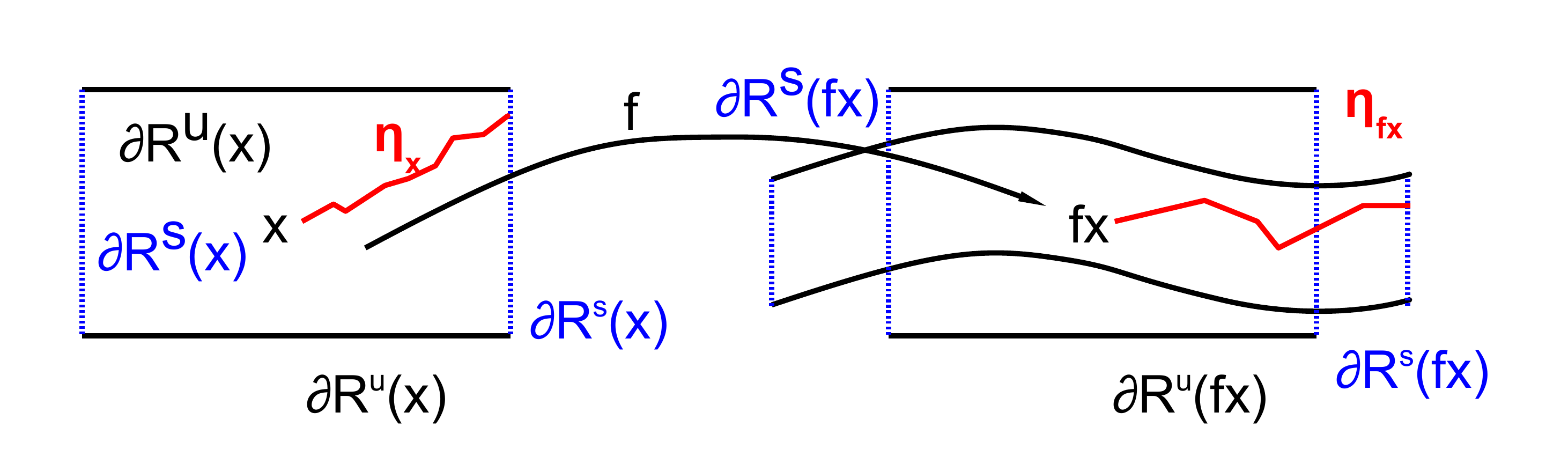}
\end{center} 

\demo
If $\eta_x(1) \in \partial  R^s(x)$, then for all $n\geq 1$, we have $\eta_{f^nx}(1)\in \partial R^s(f^nx)$. Then the map $\Ic$ defined by $\Ic(x)=1$ if $\eta_x(1)\in \partial R^s(x)$ and $\Ic(x)=0$ if not is non-decreasing along the orbits and then constant almost everywhere.

We have indeed $\int (\Ic\circ f-\Ic)d\mu=0$ and $\Ic\circ f\geq \Ic$. Hence $\Ic\circ f=\Ic $ $\mu$- a.e. and then as $\mu$ is ergodic $\Ic$ is constant $\mu$-almost everywhere. \enddemo
Assume for example that we have almost everywhere $\eta_x(1)\in \partial^s R(x)$. Hence we have $\eta_{fx}\subset f(\eta_x)$.

 The local unstable manifold at $x$ is the graph of a continuous function $g_x^u$. 

If $\eta_x=(\eta_x^1, \eta_x^2)$ we introduce the notation:
 $$\delta (x)=\max_{t\in [0, 1]} |\eta_x^2(t)-g_x^u(\eta_x^1(t))|.$$
 \begin{center}
\includegraphics[width=5.6cm]{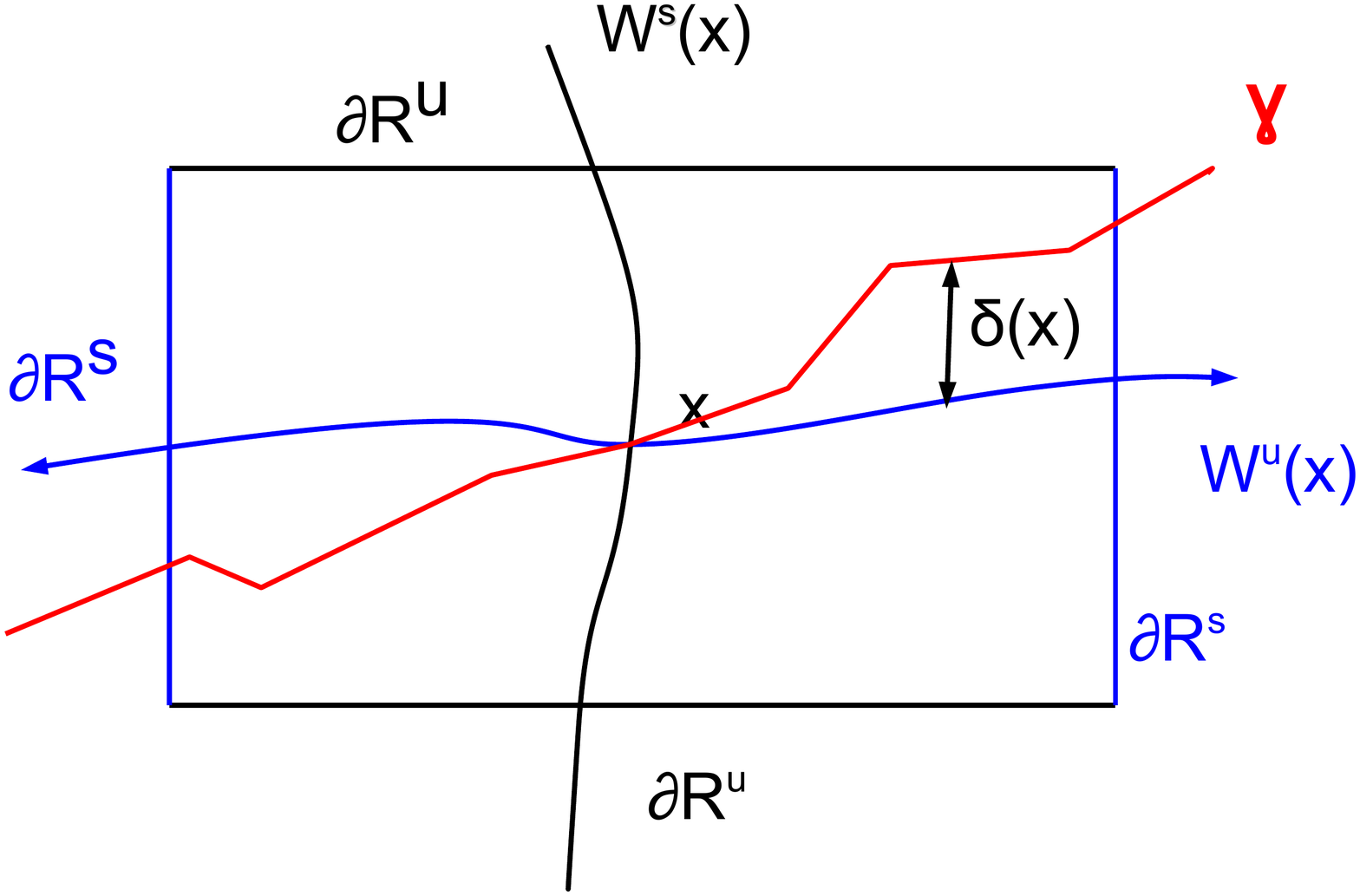}
\end{center}

Using hyperbolicity, we obtain $\delta(fx)\leq e^{-\frac{\lambda}{2}}\delta(x)$, and then $\int\delta d\mu\leq e^{-\frac{\lambda}{2}}\int\delta d\mu$ and then $\delta=0$ $\mu$ almost everywhere.

We deduce that the corresponding branch of $W^u(x)$ is contained in $\gamma$ for $\mu$-almost every $x$. 
\bigskip

Assume that $\gamma$ is irrational. Then $f_{|\gamma}$ has to be Denjoy (because for some points we have $\displaystyle{\lim_{n\rightarrow +\infty} d(f^{-n}x, f^{-n}y)=0}$).

 In this case, the only points $x\in {\rm supp}\mu$ such that $W^u(x)\not=\{ x\}$ are the endpoints of the wandering intervals and there are only  countably many such points: their set has $\mu$-measure 0. 
 
 Finally, $\gamma$ cannot be irrational. 
 
 The second part of Theorem \ref{Tirrational} is a  consequence of the following theorem that we will prove.
 
 \begin{thm}\label{thzerolyap}
Let $f:M\rightarrow M$ be a $C^1$-diffeomorphism of a manifold $M$. Let $K\subset M$ be a compact set that is invariant by $f$. We assume that $f_{|K}$ is uniquely ergodic and we denote the unique Borel invariant probability with support in $K$ by $\mu$. We assume that all the Lyapunov exponents of $\mu$ are zero. Let $x_0\in W^s(K, f)\backslash K$. Then we have:

$$ \forall \varepsilon>0, \lim_{n\rightarrow +\infty} e^{\varepsilon n}d(f^n(x_0), K)=+\infty.$$
\end{thm}
Let us now prove this theorem. 
\begin{proof}By hypothesis, we have for $\mu$-almost every point :
$$\lim_{n\rightarrow \pm\infty} \frac{1}{n}\log\| Df^n(x)\|=0.$$
We can use a refinement Kingman's subadditive ergodic theorem that is due to A.~Furman (see Theorem \ref{Tfur} of subsection \ref{ssbiLip}) that implies that we have 
$$ 
\limsup_{n\rightarrow \pm\infty} \max_{x\in K}\frac{1}{n}\log \| Df^n(x)\|\leq 0
.$$

In particular, for any $\varepsilon>0$, there exists $N\geq 1$ such that:
\begin{equation}\label{E1}\forall x\in K, \forall n\geq N, \frac{1}{n}\log\| Df^{-n}(x)\| \leq \frac{\varepsilon}{8}.\end{equation}

Observe that   the following norm with $k\ge N$ large:
\[\|u\|'_x= \sum_{n=0}^k e^{-n \varepsilon /4}\| Df^{-n}(x)u\|_x,\]
satisfies uniformly on $x$ for $u\not=0$:
\[\frac{\|D f^{-1}(x) u\|'_{f^{-1}(x)}}{\|u\|'_x}= e^{\varepsilon /4} + \frac{e^{-k \varepsilon /4}\| Df^{-k-1}(x)u\|_x-e^{\varepsilon/4}\|u\|}{\|u\|'_x}\]
\[\le e^{\varepsilon /4} + \frac{e^{-k \varepsilon /4}\| Df^{-k-1}( x)u\|_x}{\|u\|'_x}\le e^{\varepsilon /4} +e^{-k\varepsilon /8}\]
Hence by changing the Riemannian metric by the latter  one, we can assume that the norm of $D_xf^{-1}$ is smaller than $e^{ \varepsilon /3}$ for every $x\in K$.

Consequently, on a $\eta$-neighborhood $N_\eta$ of $K$, it holds for every $x\in N_\eta$ that:
\[\|D_x f^{-1}\|'\le e^{\varepsilon /2}\]

Let $x_0\in M$ be such that $x_n:=f^n(x_0)\to K$, we want to show that $$\liminf \frac1n \log d(x_n,K)\ge -\varepsilon .$$  
We suppose that $\liminf \frac1n \log d(x_n,K)< -\varepsilon $ for the sake of a contradiction. 
Hence there exists $n$ arbitrarily large so that $x_n$ belongs to the $e^{-n\varepsilon }\eta$-neighborhood of $K$.
Let $\gamma$ be a $C^1$-curve connecting $x_n$ to $K$ and of length at most $e^{-n\varepsilon }\eta$. 
By induction on $k\le n$, we notice that $f^{-k}(\gamma)$ is a curve that connects $x_{n-k}$ to $K$, and has length at most $e^{-n\varepsilon+k\varepsilon/2 }\eta$, and so is included in $N_\eta$. Thus the point $x_0$ is at most $e^{-n\varepsilon/2 }\eta$-distant from $K$. Taking $n$ large, we obtain that $x_0$ belongs to $K$. A contradiction.

\end{proof}

\enddemo

\subsection{Instability zones and the second Birkhoff theorem}\label{ssinstability}

As now we know how the dynamics restricted to $\Ic(f)$ is, we will look to the complement $\Uc(f)$ of $\Ic(f)$.

\begin{defi}{\rm 
An {\em essential curve} is a $C^0$-embedded circle in $\A$ that is not homotopic to a point, i.e. a loop that winds around the annulus.\\
An {\em  essential subannulus} of $\A$ is a subset of $\A$ that is homeomorphic to $\A$  and that contains an essential curve of $\A$.}
\end{defi}

\begin{propos}
Let $f$ be a conservative twist map. Every connected components of $\Uc(f)$ is either a bounded disc or an essential sub-annulus of $\A$.
\begin{enumerate}
\item[$\bullet$] When such a component is a disc $\Dc$ , then this disc is periodic i.e. there exists $N\geq 1$ such that $f^N(\Dc)=\Dc$. Moreover, the boundary of $\Dc$ is the union of parts of two invariant continuous graphs that have the same rational rotation number.
\item[$\bullet$]  When such a component is an essential sub-annulus, then it is invariant by $f$, and each of the two components of its boundary is either $\T\times \{\pm\infty\}$ or an invariant continuous graph.
\end{enumerate}
\end{propos}

\demo Let $U$ be a connected component of $\Uc (f)$. Then there is a partition of the set of the invariant continuous graphs in two parts: the set $\Sc_+$  of such curves that are above $U$ and the set $\Sc_-$ of those that are under $U$. Let us differentiate which cases can occur
\begin{enumerate}
\item if $\Sc_-=\Sc_+=\emptyset$, then $U=\A$ is an essential annulus;
\item if $\Sc_-=\emptyset$ and $\Sc_+\not=\emptyset$ (resp.  $\Sc_+=\emptyset$ and $\Sc_-\not=\emptyset$ ), let us denote by $\gamma_+$ (resp. $\gamma_-$) the smallest element in $\Sc_+$ (resp. the largest element in $\Sc_-$). Then $U$ is the component under $\gamma_+$ (resp. above $\gamma_-$), that is an essential sub-anulus, and its boundary is $\gamma_+$ (resp. $\gamma_-$);
\item if  $\Sc_-\not=\emptyset$ and $\Sc_+\not=\emptyset$,  let us denote by $\gamma_+$ (resp. $\gamma_-$) the smallest element in $\Sc_+$ (resp. the largest element in $\Sc_-$). Then $U$ is a connected component of the points that are between $\gamma_-$ and $\gamma_+$. If $\gamma_-\cap \gamma_+\not=\emptyset$, it is a disc $\Dc$ such that $\partial\Dc\subset \gamma_-\cup\gamma_+$; moreover, as $\gamma_-$ meets $\gamma_+$, this two curve have the same rotation number and  $\gamma_-\cap\gamma_+$ contains exactly two points of $\partial \Dc$ and they are periodic: the rotation number is rational . If $\gamma_-\cap\gamma_+=\emptyset$, then $U$ is an essential sub annulus with boundary $\gamma_-\cup\gamma_+$.
\end{enumerate}
From the fact that the invariant curves are invariant, we deduce that the the annular components of $\Uc (f)$ are invariant. The components $U$ that are homeomorphic to a disc are between two invariant curves, hence contained in an invariant domain with finite Lebesgue measure. This implies that for some $N\geq 1$, we have $f^N(U)\cap U\not=\emptyset$ and then $f^N(U)=U$.
\enddemo

\begin{defi}{\rm 
If $f$ is a conservative twist map, an annular component of $\Uc (f)$ is called an {\em instability zone}.}
\end{defi}
The following result, which was proved independently by J.~Mather (see \cite{Mat2} where the author uses variational methods) and P.~Le~Calvez (see \cite{LeC1}  where the author uses topological methods), explains why these regions are called instability zones.

\begin{thm}\label{TZI} {\bf (P.~Le~Calvez; J.~N.~Mather)}
Let $\Ac$ be an instability zone of a conservative twist map $f$ of the annulus. We choose boundaries $\Cc_-, \Cc_+$ of $\Ac$. Then there exists $x\in \Ac$ so that $\displaystyle{\lim_{k\rightarrow \pm\infty}d(f^kx, \Cc_\pm)=0}$.
\end{thm}
\begin{remas}{\rm
\begin{enumerate}
\item Note that we can choose $\Cc_-=\Cc_+$.
\item Theorem \ref{TZI} tells us that $W^u(\Cc_-)\cap W^s(\Cc_+)\cap \Ac\not=\emptyset$
\begin{center}
\includegraphics[width=2.5cm]{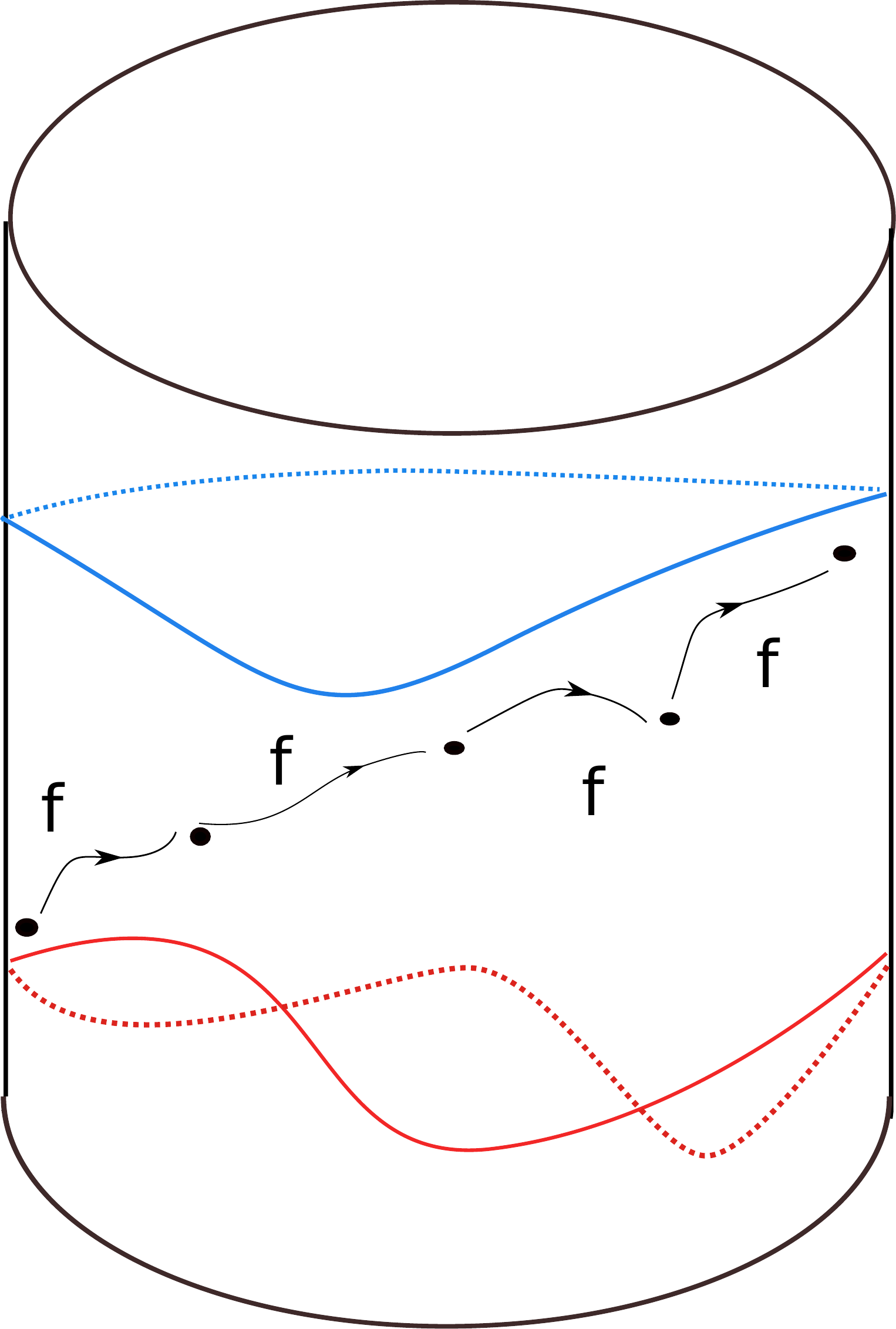}
\end{center} 
\end{enumerate}}
\end{remas}

\idemo  Let us explain in a few words what are the ideas to prove a weaker but related result due to Birkhoff: assume $\Cc_-\not=\Cc_+$,  fix a neighborhood $\Uc_-$ of  $\Cc_-$ and $\Uc_+$ of $\Cc_+$  in $\bar\Ac$,  then there exists $x\in \Uc_-$ and $N\geq 0$ so that $f^Nx\in\Uc_+$. 

The main argument is a theorem due to Birkhoff.

\begin{thm}\label{Tbir2} {\bf (G.~D.~Birkhoff)}
Let $\Ac\subset\A$ be an essential sub-annulus that is invariant by a conservative twist map of the annulus and that is equal to the interior of its closure.  Then every bounded connected component of $\partial \Ac$ is the graph of a Lipschitz map.
\end{thm}
A complete proof of Theorem \ref{Tbir2} can be found in the appendix of the first chapter of \cite{He1} (in French). \\
Then assume that $\Uc_-$ is annular and that the result we want to prove is false. For every $n\in\N$, let $V$ be the connected component of the complement in $\bar \Ac$ of $\displaystyle{\overline{\bigcup_{n\in\N}f^n(\Uc_-)}}$ that contains $\Cc_+$. One can check that the interior of $\bar V$ satisfies the hypothesis of Theorem \ref{Tbir2}, hence we find an invariant continuous graph that is in $\Ac$ (the boundary of $V$), that is incompatible with the definition of an instability zone.\enddemo

Note an important corollary of theorem \ref{Tbir2}.

\begin{corol}
Let $\gamma $ be an essential curve that is invariant by a conservative twist map. Then $\gamma$ is the graph of a Lipschitz map.
\end{corol}

\begin{exam}\label{ExZI} {\rm This example was introduced by Birkhoff in \cite{Bir2}. We consider the Hamiltonian flow $f$ of the pendulum for a small enough time. Using a perturbation of the generating function of $f$, we can create a transverse intersection between the lower stable branch and the lower unstable branch of the hyperbolic fixed point: 
\begin{center}
\includegraphics[width=3cm]{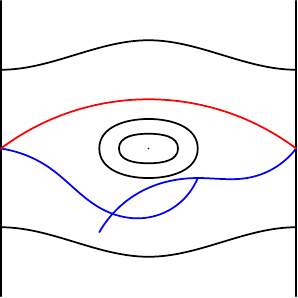}
\end{center}
Then the remaining separatrix is the upper boundary of an instability zone.}
\end{exam}
\begin{exer}{\rm 
Prove the last assertion in Example \ref{ExZI}.}
\end{exer}
Michel Herman proved in \cite{He2} that for a general conservative twist map, there is no essential invariant curve that contains a periodic point. More precisely:\\
{\em Let $k\in [1, +\infty]$ be a positive integer or $\infty$. There exists a dense $G_\delta$-subset $\Gc$ of the set of the $C^k$ PSTM such that every $f\in\Gc$ has no invariant essential curve that contains a periodic point.}\\
The proof of this result is proposed in Exercice \ref{Exo}.

\begin{quest} {\rm For which parameters $\lambda$ does the standard map $f_\lambda$ satisfy this property? }\end{quest}

\begin{quest} {\rm How is a ``general'' boundary of an instability zone? Is it the boundary of one or two intability zone(s)? Is it smooth? How is its rotation number: Diophantine, Liouville? }\end{quest}

\begin{rema}{\rm This result of Michel Herman joined to the fact that there exist open sets of $C^\infty$ conservative twist maps that have a lot of  (Diophantine) invariant graphs, allows us to state :} \end{rema}

\begin{propos}\label{Pstableirrat}
There exists a dense $G_\delta$-subset $\Gc$ (for the $C^\infty$-topology) in a non-empty open set of conservative $C^\infty$ twist map such that every $f\in\Gc$ has a bounded instability zone with irrational boundaries.
\end{propos}

Then the stable set of such an irrational boundary is not empty (because of Theorem \ref{TZI}) but the convergence to such a boundary is slower than exponential (because of Theorem \ref{Tirrational}).

\begin{exer}{\rm 
Prove Proposition \ref{Pstableirrat}.}
\end{exer}

\begin{quest}{\rm 
 For which parameters $\lambda$ has the standard map $f_\lambda$ an irrational boundary of instability zone?}

\end{quest}

%Parler de la nong\'eg\'erivit\'e des courbes irrationnelles?(parler alors points p�riodiques g�n�riques)

%Existence d'ouverts ayant des zones d'instabilit�...
%le theorem varit� stable, instable

%les exemples: moi (construction � la herman).

%question: 2 cotes zone instabilit�???

%nb rota, irregularite?

\section{Aubry-Mather theory}\label{SAubryMather}

\subsection{ Action functional and minimizing orbits}\label{ssminimizing}
In this section, we assume that $S:\R^2\rightarrow \R$ is a generating function of a lift $F:\R^2\rightarrow \R^2$ of a conservative twist map $f:\A\rightarrow \A$.

\begin{defi}{\rm 
If $k\geq 1$, one defines the {\em action functional} $\Fc_{k+1}:\R^{k+1}\rightarrow \R$ by $$\displaystyle{\Fc(\theta_0, \dots , \theta_k)=\sum_{j=1}^kS(\theta_{j-1}, \theta_j)}.$$}
\end{defi}

For every $k\geq 2$ and every $\theta_b, \theta_e\in \R^n$, the function $\Fc_{k+1}$ (or $\Fc$) restricted to the set $\Ec (k+1, \theta_b, \theta_e)$ of $(k+1)$-uples $(\theta_0, \dots ,\theta_k)$ beginning at $\theta_b$ and ending at $\theta_e$, i.e. such that $\theta_0=\theta_e$ and $\theta_k=\theta_e$, has a minimimum and at every critical point for $\Fc_{k+1|\Ec (k+1, \theta_b, \theta_e)}$, the following sequence is a piece of orbit for $F$:
 $$(\theta_0, -\frac{\partial S}{\partial \theta}(\theta_0, \theta_1)), (\theta_1,  \frac{\partial S}{\partial \Theta}(\theta_0, \theta_1)), (\theta_2,    \frac{\partial S}{\partial \Theta}(\theta_1, \theta_2)), \dots , (\theta_k,  \frac{\partial S}{\partial \Theta}(\theta_{k-1}, \theta_k)).$$
 Observe that for such a critical point, we have $\frac{\partial S}{\partial \Theta}(\theta_{i-1}, \theta_i)+\frac{\partial S}{\partial \theta}(\theta_{i}, \theta_{i+1})=0$ for every $0<i<k$.
 
 \begin{exam}{\rm 
 To illustrate the notion of generating function, let us introduce a very classical example of twist map that is due to G.D.~Birkhoff: the so-called Birkhoff billiard. Play billiard on a planar billiard table with a $C^2$ and convex boundary with non-vanishing curvature. Then we can choose symplectic coordinates (angular coordinate for the point of bounce and radial coordinate that is the sinus of the angle of reflection) in such a way that the dynamical system becomes a conservative twist map (see \cite{Sib} for details). 
 \begin{center}
\includegraphics[width=3cm]{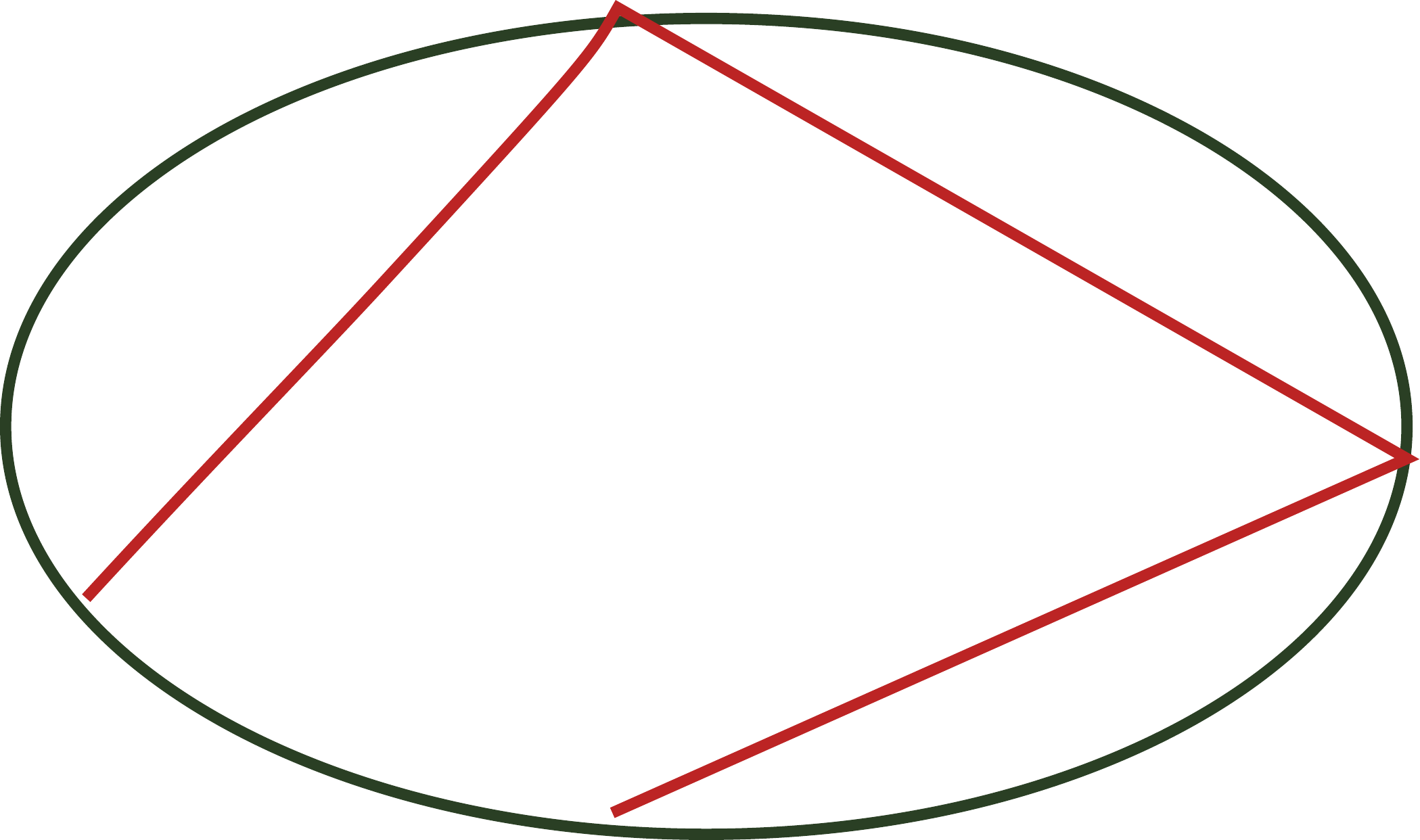}
\end{center}

 In these coordinates, if $\theta_0, \dots, \theta_n\in \R^{n+1}$, then $\Fc(\theta_0, \dots, \theta_n)$ is just the length of the polygonal line that joins the successive points with angular coordinates $\theta_0, \dots, \theta_n$.

 }
 
  \end{exam}
 
 \begin{defi} {\rm A finite or infinite sequence of real numbers $(\theta_n)_{n\in J} $ is a {\em minimizer} if for every segment $[\ell , k]\subset J$, $(\theta_n)_{\ell\leq n\leq k}$ is  a global minimizer of  $\Fc_{k-\ell+1|\Ec (k-\ell+1, \theta_\ell, \theta_k)}$.
 
 When $J=\Z$, we say that $(\theta_n)$ is a {\em minimizing sequence}; we denote the set of minimizing sequences by $\Mc\subset \R^\Z$.}

An orbit $(\theta_n, r_n)$ of $F$ (and by extension its projection on $\A$) is {\em minimizing} if its projection $(\theta_n)$ is a minimizing sequence. \\
   \end{defi}

\begin{rema} {\rm Observe that a minimizer is always the projection of a piece of orbit. 
From Lemma \ref{Lcoercive}, we can deduce
\begin{enumerate}
\item[$\bullet$] in every $\Ec=\Ec (k+1, \theta_b, \theta_e)$, there exists a minimizer of $\Fc_{|\Ec}$; such a minimizer is a segment of  the projection of an (non necessarily minimizing) orbit;
\item[$\bullet$] if $(q, p)\in \Z^*\times\Z$, the restriction of $\Fc_{q+1}$ to the set $\{ (\theta_k);\theta_{k+q}=\theta_k+p\}$ has a global minimizer. Any such minimizer is the projection of an orbit and we will even see in Proposition \ref{Ppermin} that it is a minimizing sequence.
\end{enumerate}}
\end{rema}

 The following theorem is due to J.~Mather and proved in subsection \ref{ssgraphmin}.
 
 \begin{thm} {\bf (J.~N.~Mather)}
 Assume that the graph of a continuous map $\psi:\T\rightarrow \R$ is invariant by a conservative twist map $f$. Then for any generating function associated to $f$, all the orbits contained in the graph of $\psi$ are minimizing.
 \end{thm}

% \subsection{Comparison of different minimizers and some consequences}

Now we will give some properties of the minimizers and prove the  existence of some periodic minimizers.
 
 \begin{propos}\label{noncrossing} {\bf (Aubry \& Le Daeron non-crossing lemma)} Assume $(b-a)(B-A)\leq 0$. Then 
 $$S(a, A)+S(b, B)-S(a, B)-S(b, A)\geq 0$$
 and equality occurs if and only if  $(b-a)(B-A)= 0$.
 \end{propos}
 
 \begin{center}
\includegraphics[width=5cm]{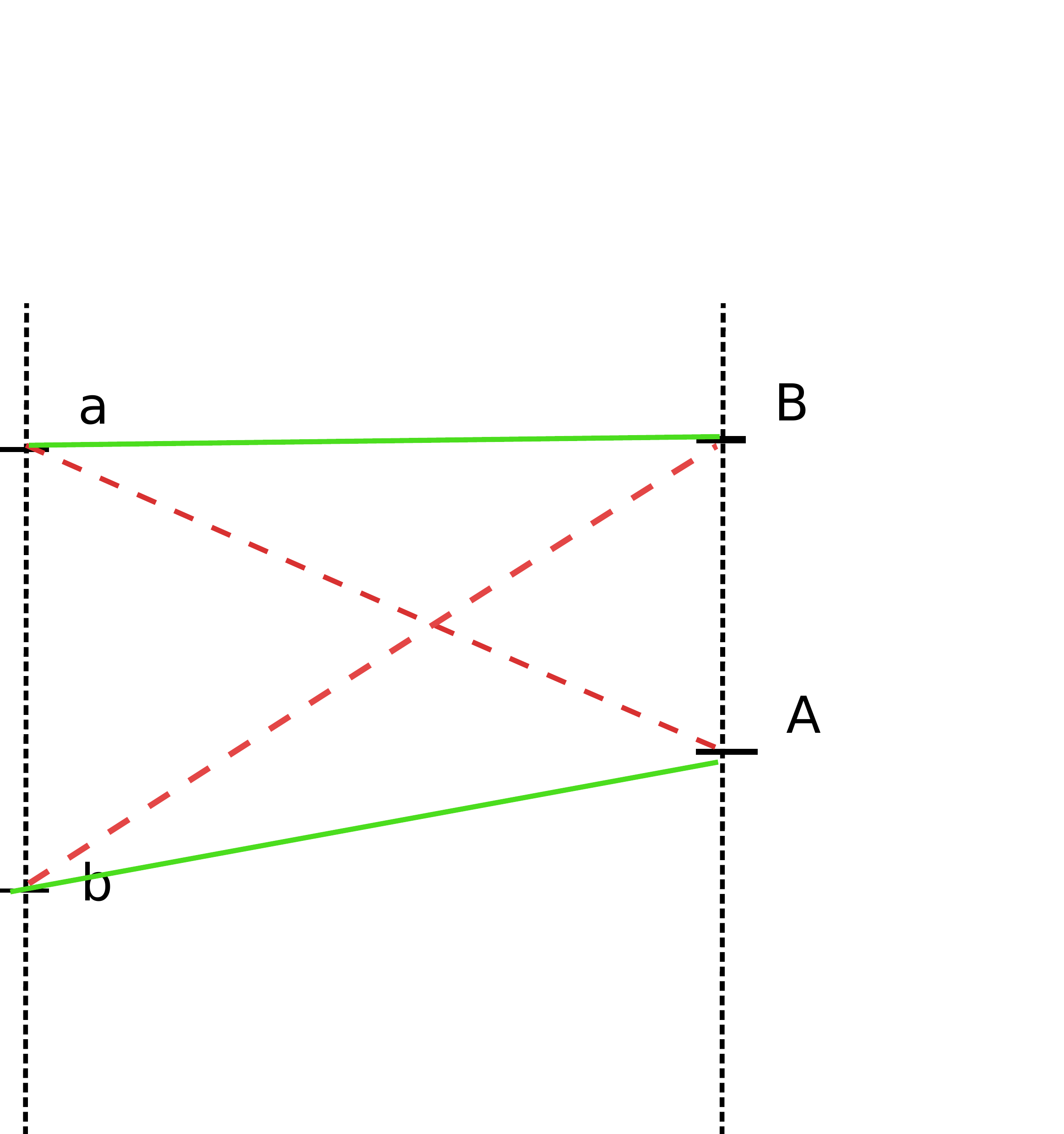}
\end{center}

\demo
Let us use the notation $A_t=A+t(B-A)$ and $a_t=a+t(b-a)$. We have: 
$$\begin{matrix}S(a, A)+S(b, B)-S(a, B)-S(b, A)&=\left( S(b, B)-S(b, A)\right)-\left( S(a, B)-S(a, A)\right)\hfill\\
&=(B-A) \int_0^1\left( \frac{\partial S}{\partial \Theta} (b, A_t)-\frac{\partial S}{\partial \Theta}(a, A_t)\right)dt\hfill\\
&=(b-a)(B-A)\int_0^1\int_0^1\frac{\partial^2 S}{\partial \theta\partial \Theta}(a_s, A_t)ds.dt.\hfill
\end{matrix}
$$
From $\frac{\partial^2 S}{\partial \theta\partial \Theta}<0$, we deduce the wanted result.
\enddemo
\begin{defi}{\rm
If $(\theta_k)$ is a finite or infinite sequence of real numbers, its {\em Aubry diagram} is the graph of the function obtained when interpolating linearly the sequence $(k, \theta_k)$.

Two sequences $(a_k)_{k\in I}$ and $(b_k)_{k\in I}$ {\em cross} if for some $k$, $j$: $(a_k-b_k)(a_{j}-b_{j})<0$.}
\end{defi}

\begin{rema}{\rm They are two types of crossing: at an integer or at a non-integer:
 \begin{center}
\includegraphics[width=7cm]{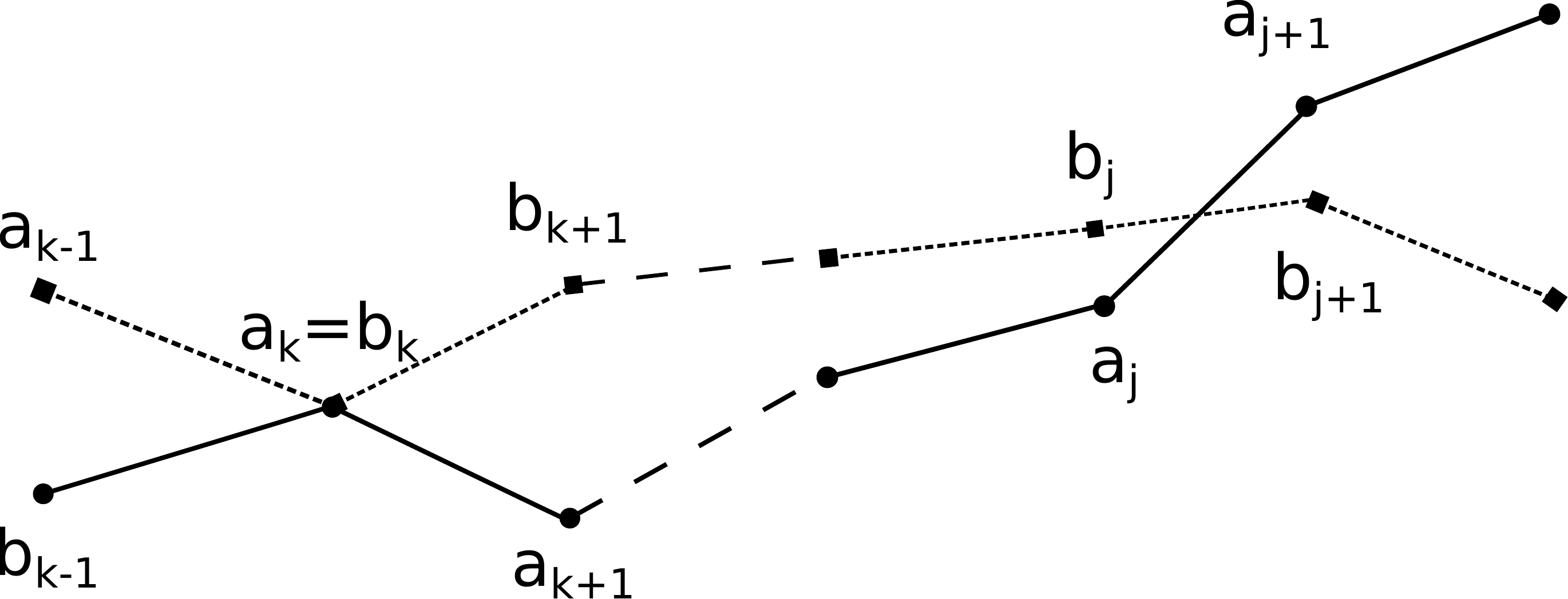}
\end{center}
Note that if two distinct minimizers are such that for a $k$ we have $a_k=b_k$, then we have $a_{k-1}\not= b_{k-1}$ and $a_{k+1}\not=b_{k+1}$; indeed, if two successive terms coincide, then they correspond to a same orbit and then to the same minimizer.}
\end{rema}

\begin{propos}\label{PAubryfund}{\bf (Aubry fundamental lemma)} Two distinct minimizers cross at most once.

\end{propos}
\demo
Assume that the minimizers $(a_k)$ and $(b_k)$ cross at two different times $t_1$ and $t_2$.  Let us introduce the notation $k_i=[t_i]$. We consider the the following finite segments:
\begin{enumerate}
\item[$\bullet$] $A=(a_k)_{k_1\leq k\leq k_2+1}$;
\item[$\bullet$] $B=(b_k)_{k_1\leq k\leq k_2+1}$;
\item[$\bullet$] $\alpha=(a_{k_1}, b_{k_1+1}, \dots, b_{k_2}, a_{k_2+1})$; 
\item[$\bullet$] $\beta=(b_{k_1}, a_{k_1+1}, \dots, a_{k_2}, b_{k_2+1})$.
\end{enumerate}
If $t_1$ or $t_2$ is not an integer, we deduce from Proposition \ref{noncrossing} that 
$$\begin{matrix}\Fc(A)+\Fc(B)- \Fc(\alpha)-\Fc(\beta)=\hfill{}\quad\\
\quad{}\hfill\displaystyle{\sum_{i=1}^2\left( S(a_{k_i}, a_{k_i+1})+S(b_{k_i}, b_{k_i+1})-S(a_{k_i}, b_{k_i+1})-S(b_{k_i}, a_{k_i+1})\right)}>0.\end{matrix}$$
 As $A$ and $\alpha$ (resp. $B$ and $\beta$) have same endpoints, we deduce that $A$ or $B$ is not minimizing, and this is a contradicton.

If $t_i=k_i$ are both integers, then we obtain $\Fc(A)+\Fc(B)- \Fc(\alpha)-\Fc(\beta)=0$. As $\Fc(A)\leq \Fc(\alpha)$ and $\Fc (B)\leq \Fc(\beta)$, we deduce that $\alpha$ and $\beta$ are also minimizers. But $\alpha$ and $A$ coincides for integers $k_2$ and $k_2+1$, hence $\alpha=A$ and then $A=B$.
\enddemo
\begin{defi}{\rm 
If $(q, p)\in\N^*\times \Z$, a sequence $(\theta_n)_{n\in\Z}$ is a {\em $(q, p)$-minimizer} if
\begin{enumerate}
\item $\forall n, \theta_{n+q}=\theta_n+p$;
\item $(\theta_n)_{0\leq n\leq q-1}$ is a minimizer of the function $\displaystyle{(\alpha_n)_{0\leq n\leq q-1}\mapsto \sum_{n=0}^q S(\alpha_n, \alpha_{n+1})}$ (with the convention $\alpha_q=\alpha_0+p$).
\end{enumerate}}
\end{defi}
Observe that $(q, p)$-minimizer is the projection of an orbit $(\theta_n, r_n)$ for $F$ such that $(\theta_{n+q}, r_{n+q})=(\theta_n, r_n)+(p, 0)$. Hence it corresponds to a $q$-periodic orbit for $f$.

\begin{propos}\label{Ppermin}
Any $(q, p)$-minimizer is a minimizing sequence.
\end{propos}

\begin{exer}{\rm The goal of the exercise is to prove Proposition \ref{Ppermin}.\\
(a) Using Proposition \ref{PAubryfund}, prove that for every  $(q, p)\in\N^*\times \Z$ and $k\geq 1$, two distinct $(q, p)$-minimizers cannot cross.\\
{\em Hint: {\em prove that if they cross, they cross two times within a period.}}\\
(b) Deduce that for every  $(q, p)\in\N^*\times \Z$ and $k\geq 1$, every $(kq, kp)$-minimizer is in fact a $(q, p)$-minimizer.\\
(c) Deduce that being a $(q,p)$-minimizer is equivalent to be a $(kq, kp)$-minimizer. \\
(d) Deduce Proposition \ref{Ppermin}.}

\end{exer}

\begin{notat}{\rm If $(q,p)\in \Z^2$, we denote by $T_{q, p}:\R^\Z\rightarrow \R^\Z$ the map defined by $T_{q, p}((x_k)_{k\in\Z})=(x_{k-q}+p)_{k\in\Z}$.}\end{notat}

Note that if $(\theta_k)_{k\in\Z}$ is a $(q, p)$ minimizer, then $T_{q, p}\left( (\theta_k)_{k\in\Z}\right) =(\theta_k)_{k\in\Z}$.

\begin{corol}\label{Ccomparper}
If $(\theta_k)_{k\in\Z}$ and $(\alpha_k)_{k\in\Z}$ are two $(q, p)$-minimizers, then they don't cross. In particular, $(\theta_k)_{k\in\Z}$ and $T_{a, b}\left( (\theta_k)_{k\in\Z}\right)$ do not cross.
\end{corol}
\begin{propos}\label{Pminper}
For every $q\in \N^*$, $p\in\Z$, there exists at least one $(q, p)$-minimizer.
\end{propos} 
\demo
We assume that $S$ is a generating function of a lift $F$ of the conservative twist map $f$. 
\begin{lemm}\label{Lcoercive}
 We have $\displaystyle{\lim_{|\Theta-\theta|\rightarrow +\infty}\frac{S(\theta, \Theta)}{|\Theta-\theta|}=+\infty
 }$.
 \end{lemm}
 \demo Using the notation $\theta_t=\theta+t(\Theta-\theta)$, we have $$\begin{matrix} S(\theta, \Theta)&=S(\theta, \theta)+\int_0^1\frac{\partial S}{\partial \Theta}(\theta, \theta_t)(\Theta-\theta)dt\hfill\\
 &=S(\theta, \theta)+\int_0^1\frac{\partial S}{\partial \Theta}(\theta_t, \theta_t)(\Theta-\theta)dt-\int_0^1\int_0^t\frac{\partial^2 S}{\partial \theta\partial \Theta}(\theta_{s}, \theta_t)(\Theta-\theta)^2dsdt\\
 &\geq m-M|\Theta-\theta|+\frac{\varepsilon}{2}(\Theta-\theta)^2\hfill 
 \end{matrix}$$
 where $\displaystyle{m=\min_{\theta\in[0, 1]}S(\theta, \theta)}$ and $\displaystyle{M=\max_{\theta\in[0, 1]}\left| \frac{\partial S}{\partial \Theta}(\theta, \theta)\right|}$.
 \enddemo
 We know consider the set
 $$\Ec(q, p)=\{ (\theta_k)_{k\in\Z}; \forall k\in\Z, \theta_{k+q}=\theta_k+p\}$$
 and define $\Wc:\Ec(q,p)\rightarrow \R$ by
 $$\Wc((\theta_k)_{k\in\Z})=\sum_{k=0}^{q-1}S(x_k, x_{k+1}).$$
 Note that if $\ell\in\Z$, then $\Wc((\theta_k)_{k\in\Z})=\Wc((\theta_k+\ell)_{k\in\Z})$. Hence we can define $\Wc$ on the quotient of $\Ec(q, p)$ by the diagonal action of $\Z$. On this space, $\Wc$ is  coercive and has then a global minimimum. Then this global minimum is attained at a $(q, p)$-minimizer.

\enddemo

\begin{exer}{\rm 
Write the details in the proof of   Proposition \ref{Pminper}. }
\end{exer}

\subsection{$F$-ordered sets}
\begin{defi}{\rm 
We say that a subset $E\subset \R^2$ is $F$-ordered if it is invariant by $F$ and every integer translations $(\theta, r)\mapsto (\theta+k, r)$ with $k\in\Z$ and if
$$\forall x, x'\in E, \pi(x)<\pi(x')\Rightarrow \pi\circ F(x)<\pi\circ F(x').$$}
\end{defi}

\begin{rema}{\rm We deduce from Corollary \ref{Ccomparper} that if $q\in\Z^*$ and $p\in\Z$, the union of the $(q,p)$-minimizing orbits is an $F$-ordered set.}\end{rema}

\begin{exer}{\rm
Let $\psi:\T\rightarrow \R$ be a continuous map such that the graph of $\psi$ is invariant by a conservative twist map $f$. Prove for any lift $F$ of $f$, the graph of $\psi$ is $F$-ordered.}
\end{exer}

The following proposition explains how we can construct other $F$-ordered sets.

\begin{propos}\label{PForderedlimit}
Let $F$ be a lift of a conservative twist map.
\begin{enumerate}
\item The closure of every $F$-ordered set is $F$-ordered;
\item Let $(E_n)_{n\in\N}$ be a sequence of $F$-ordered sets. Let $E\in\R^2$ be the set of points $x\in\R^2$ so that there exist $(x_n)\in\R^2$ satisfying $x_n\in E_n$ and $\displaystyle{\lim_{n\rightarrow \infty} x_n=x}$.
\end{enumerate}
Then $E$ is $F$-ordered.
\end{propos}
\begin{rema}{\rm  The main remark that is useful to prove Proposition \ref{PForderedlimit} is the following one.\\
Assume that $E\subset \R^2$ is invariant by $F$ and all maps $(\theta, r)\mapsto (\theta+k, r)$ with $k\in\Z$. Then $E$ is $F$-ordered if and only if 
$$\forall x, x'\in E, \pi(x)<\pi(x')\Rightarrow \pi\circ F(x)\leq\pi\circ F(x')\quad{\rm and}\quad \pi\circ F^2(x)\leq\pi\circ F^2(x').$$
To prove that, observe that if $\pi\circ F(x)=\pi\circ F(x')$ for some $x\not=x'$ in $\R^2$, then $(\pi\circ F^{-1} (x), \pi\circ F^{-1}(x'))$ and $(\pi\circ F (x), \pi\circ F(x'))$ are not in the same order.}\end{rema}
\begin{propos}\label{PAMgraph}
Let $F$ be a lift of a conservative twist map and let $E\subset \R^2$ be a non-empty and closed $F$-ordered set. Then $\pi$ maps $E$ homeomorphically onto a closed  subset of $\R$ that is invariant by the map $t\in\R\mapsto t+1$.
\end{propos}
\demo
The map $\pi$ is continuous and open. Assume that there exist two points  $x\not=y$ of $E$ such that $\pi(x)=\pi(y)$. Because of the twist condition, we have $x_-=\pi\circ F^{-1}(x)\not=\pi\circ F^{-1}(y)=y_-$ and this contradicts the fact that $E$ is $F$-ordered.

We just have to prove that $\pi(E)$ is closed. Assume that $(x_n)$ is a sequence of points of $E$ such that $(\pi(x_n))$ converges to some $\theta\in \R$. Then there exists $a, b\in\Z$ so that $\forall n\in \N, \pi(x_0)+a<\pi(x_n)<\pi(x_0)+b$. Because $E$ is $F$-ordered, we have then $\forall n\in \N, \pi\circ F(x_0)+a<\pi\circ F (x_n)<\pi\circ F(x_0)+b$. Hence 
$$x_n\in \pi^{-1}([\pi(x_0)+a, \pi(x_0)+b])\cap F^{-1}(\pi^{-1}([\pi\circ F(x_0)+a, \pi\circ F(x_0)+b]))=K.$$
 \begin{center}
\includegraphics[width=5cm]{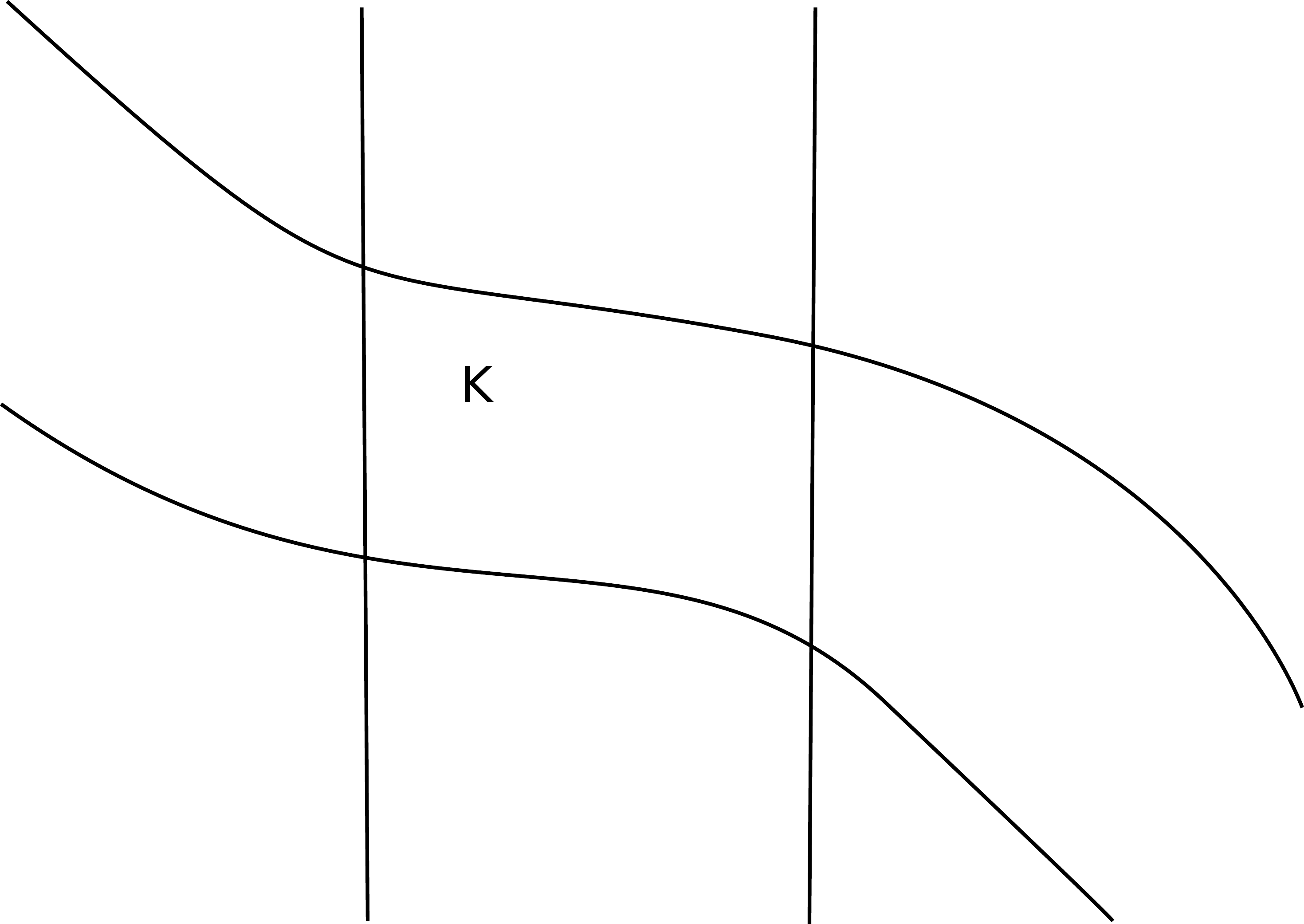}
\end{center}

Because of the twist condition, $K$ is compact. Hence we can extract a convergent subsequence from $(x_n)$. Because $E$ is closed, $x=\lim x_n\in E$ and then $\theta=\pi(x)\in \pi(E)$.
\enddemo
We deduce the following statement.

\begin{propos} \label{PAMhomeo}
Let $F$ be the lift of a conservative twist map and let $E\subset \R^2$ be a non-empty and closed $F$-ordered set. Then there exists an increasing   homeomophism $H:\R\rightarrow \R$ such that
\begin{enumerate}
\item[$\bullet$] $\forall t\in \R, H(t+1)=H(t)+1$;
\item [$\bullet$] $\forall x\in E, H\circ \pi(x)=\pi\circ F(x)$.
\end{enumerate}
\end{propos}
Hence the dynamics $F$ restricted to $E$ is conjugated (via $\pi$) to the one of a lift of a circle homeomorphism. We even deduce from Proposition \ref{PAMLip} that $H$ is bi-Lipschitz. We can then associate to every $F$-ordered set a rotation number.

\begin{propos}\label{PAMLip}
Let $f:\A \rightarrow \A$ be a conservative twist map and $x\in \A$. Then there exists a $C^1$-neighborhood $\Uc$ of $f$, a neighborhood $U$ of $x$ in $\A$ and a constant $C>0$ such that 

if $E\subset\R^2$ is a $G$-ordered set for a lift $G$ of some $g\in\Uc$ that meets $U+\Z\times\{ 0\}$, then $E$ is the graph of some $C$-Lipschitz map $\psi:\pi(E)\rightarrow \R$.

\end{propos}
Note that this proposition is similar to Theorem \ref{Tbir1} (in fact, we can  deduce Theorem \ref{Tbir1} from Proposition \ref{PAMLip}). \demo Let $F$ be a lift of the conservative twist map $f=(f_1, f_2)$, let $\varepsilon>0$ be so that $\frac{\partial f_1}{\partial r}\in (\varepsilon, \frac{1}{\varepsilon})$ and let $x=(\theta, r)$ be a point of $\R^2$. Let us choose a compact neighbourhood $B$ of $x$. \\
Then for every $y=(\alpha, \rho)\in B$, if we use the notation $y_-=F^{-1}(y)=(\alpha_-, \rho_-)$ and $y_+=F(y)=(\alpha_+, \rho_+)$,  the curves   $F^{-1}(\{\alpha_+\}\times[r_+-\frac{1}{\varepsilon}, r_++\frac{1}{\varepsilon}])$ and  $F(\{\alpha_-\}\times [ r_--\frac{1}{\varepsilon},r_-+\frac{1}{\varepsilon}])$ are graphs of some $C^1$ functions $v_{y,-}, v_{y, +}$ whose domains contain $[\alpha-1, \alpha+1]$. 
 \begin{center}
\includegraphics[width=7cm]{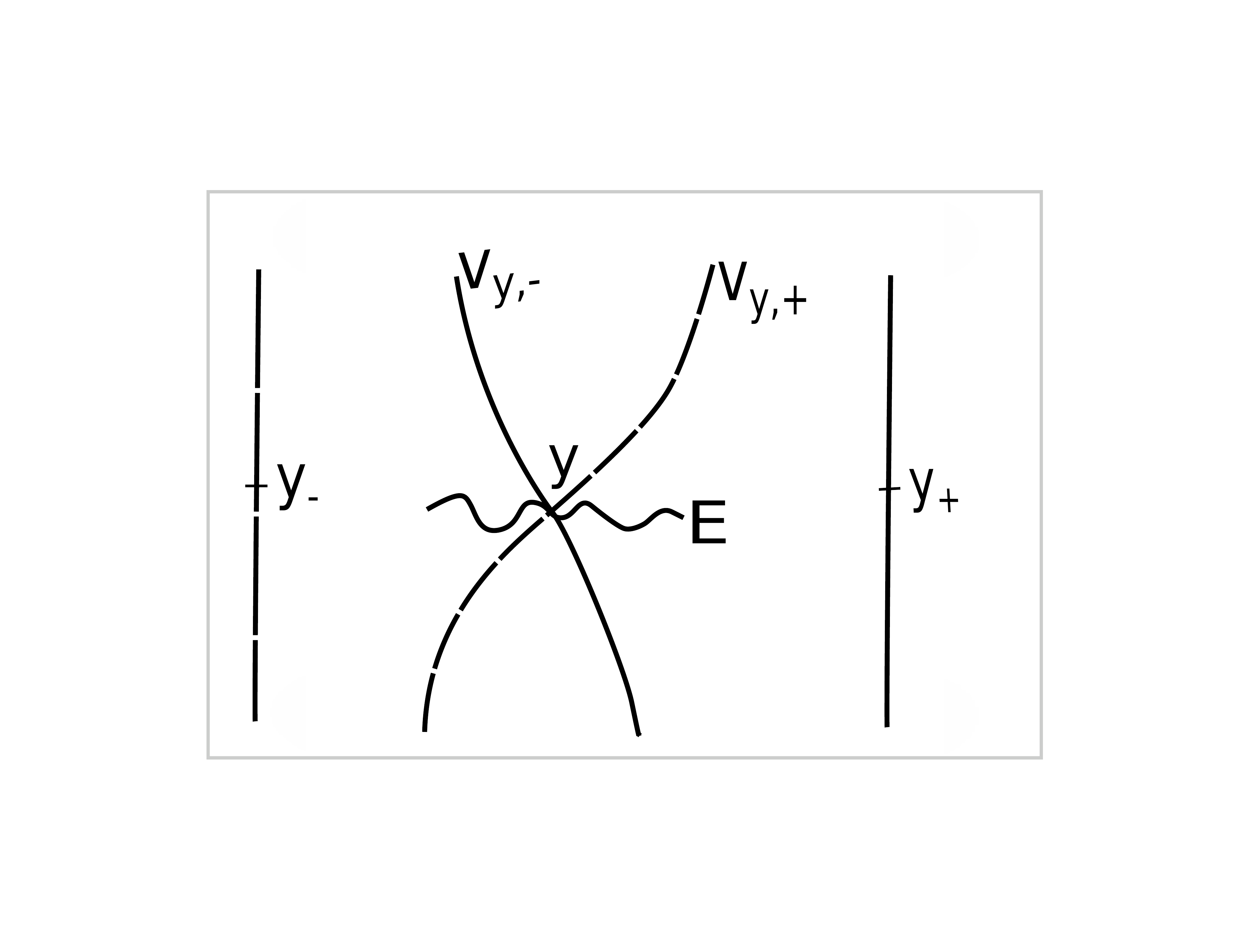}
\end{center}

Because $F(F^{-1}(B)+\{0\}\times[-\frac{1}{\varepsilon}, \frac{1}{\varepsilon}])$ and $F^{-1}(F(B)+\{0\}\times[-\frac{1}{\varepsilon}, \frac{1}{\varepsilon}])$ are compact, there exists $K>0$ such that $\R\times [-K, K]$ contains these two sets.

We define now $\Uc$ as being the set of conservative twist maps $g=(g_1, g_2)$ with a lift $G$ such that
\begin{enumerate}
\item[$\bullet$]  $\forall x \in \left( G^{-1}(B)+\{0\}\times[-\frac{1}{\varepsilon}, \frac{1}{\varepsilon}]\right) \cup G^{-1}\left( G(B)+\{0\}\times[-\frac{1}{\varepsilon}, \frac{1}{\varepsilon}]\right), \frac{\partial g_1}{\partial r}(x)\in (\varepsilon, \frac{1}{\varepsilon})$;
\item[$\bullet$]  $G( G^{-1}(B)+\{0\}\times[-\frac{1}{\varepsilon}, \frac{1}{\varepsilon}]) \cup G^{-1}( G(B)+\{0\}\times[-\frac{1}{\varepsilon}, \frac{1}{\varepsilon}])\subset \R\times [-K, K]$.
\end{enumerate}
Assume that $G$ is such a lift of $g\in\Uc$. Let $E$ be a $G$-ordered set that meets $B$ at some $y$. We deduce from Proposition \ref{PAMgraph} that $E$ is the graph of a map $\psi:\pi(E)\rightarrow \R$ and then $y=(\alpha, \psi(\alpha))$ for some $\alpha\in\pi(E)\subset \R$. Because $g\in\Uc$, we know that  $G( G^{-1}(y)+\{0\}\times[-\frac{1}{\varepsilon}, \frac{1}{\varepsilon}])$ and $G^{-1}( G(y)+\{0\}\times[-\frac{1}{\varepsilon}, \frac{1}{\varepsilon}])$ are some subsets of $\R\times [-K, K]$ and are graphs of some $C^1$ maps $v_-$, $v_+$ whose domains contain $[\alpha -1, \alpha+1]$. We can even extend these functions to $\R$ by asking that $v_-$ (resp. $v_+$) is the graph of $G^{-1}(\Vc(G(y)))$ (resp. $G(\Vc(G^{-1}(y)))$).\\
Because $E$ is $G$-ordered, we have $G^{-1}(\{z\in  E, \pi(z)\leq \pi\circ G(y)\})=\{ z\in E; \pi(z)\leq \alpha\}$. Hence $\{ z\in E; \pi(z)< \alpha\}$ is in the connected component of $\R^2\backslash G^{-1}(\Vc(G(y)))$ that is under $v_-$. Using some similar arguments, we finally obtain 
$$\forall t\in (-\infty, \alpha)\cap \pi(E), v_+(t)<\psi(t)<v_-(t)$$
 and $$\forall t\in (\alpha, +\infty)\cap \pi(E), v_-(t)<\psi(t)<v_+(t).$$
 Using the invariance by integer translation of $E$ (i.e. $E+(1, 0)=E$) and the fact that the graphs of $v_-$ and $v_+$ restricted to $[\alpha-1, \alpha+1]$ are contained in $\R\times [-K, K]$, we deduce that $E\subset \R\times [-K, K]$.
 
 We will now add a condition to define $\Uc$. Let $L>\frac{1}{\varepsilon}$ be a real number such that
 $$\forall x \in F^{-1}(\R\times [-K-\frac{1}{\varepsilon}, K+\frac{1}{\varepsilon}])\cup  (\R\times [-K-\frac{1}{\varepsilon}, K+\frac{1}{\varepsilon}]), \max \{\left|\frac{\partial f_2}{\partial r}(x)\right|, \left|\frac{\partial f_1}{\partial \theta}(x)\right|\}<L.$$
 
 Then we ask that every lift $G$ of an element $g=(g_1, g_2)$ of $\Uc$ (in addition to the other conditions we gave before that) satisfies
 \begin{enumerate}
\item[$\bullet$] we have  $$\forall x \in G^{-1}(\R\times [-K-\frac{1}{\varepsilon}, K+\frac{1}{\varepsilon}])\cup (\R\times [-K-\frac{1}{\varepsilon}, K+\frac{1}{\varepsilon}]), \max \{\left|\frac{\partial g_2}{\partial r}(x)\right|, \left|\frac{\partial g_1}{\partial \theta}(x)\right|\}<L;$$
\item[$\bullet$] and $$\forall x \in G^{-1}(\R\times [-K-\frac{1}{\varepsilon}, K+\frac{1}{\varepsilon}])\cup(\R\times [-K-\frac{1}{\varepsilon}, K+\frac{1}{\varepsilon}]), \frac{\partial g_1}{\partial r}(x)>\varepsilon.$$
 \end{enumerate}

Let us now consider $y=(\alpha, \psi(\alpha))\in E$. Repeating the same argument than before, we know that 
$$\forall t\in \pi(E), \min\{ v_-(t), v_+(t)\}\leq \psi(t)\leq  \max\{ v_-(t), v_+(t)\}.$$
 Note that $v_-'(t)=-\frac{\partial g_1}{\partial\theta}(t, v_-(t))\left( \frac{\partial g_1}{\partial r}(t, v_-(t))\right)^{-1}$ and then for every $t\in [\alpha-1, \alpha+1]$, $|v'_-(t)|<\frac{L}{\varepsilon}$. \\
 Moreover, we have $v_+'(t)=\frac{\partial g_2}{\partial r}(G^{-1}(t, v_+(t)))\left( \frac{\partial g_1}{\partial r}(G^{-1}(t, v_+(t)))\right)^{-1}$ and then for every $t\in [\alpha-1, \alpha+1]$, $|v'_+(t)|<\frac{L}{\varepsilon}$.\\
  We introduce the notation $C=\frac{L}{\varepsilon}$. We have then:
 $\forall t\in [\alpha-1, \alpha+1], |\psi (t)-\psi(\alpha)|\leq \max\{ |v_-(t)-\psi(\alpha)|, |v_+(t)-\psi(\alpha)|\}\leq C|t-\alpha|$.\\
 This implies that $\psi$ is $C$-Lipschitz.

\enddemo

 \subsection{ Aubry-Mather sets}
 
 \begin{defi}{\rm
 Let $F$ be a lift of a conservative twist map $f$. An {\em Aubry Mather set} for $F$ is a closed $F$-ordered set.\\
 The Aubry-Mather set is {\em minimizing} if every orbit contained in it is minimizing.}
 \end{defi}
 
 We noticed that any $F$-ordered set has a rotation number.  
 
 \begin{notat}{\rm  If $E$ is an Aubry-Mather set, we denote by $\rho(E)$ its rotation number.
 The Aubry-Mather set $E$ is said to be {\em rational} (resp. {\em irrational}) if $\rho(E)$ is rational (resp. irrational).}\end{notat}
 
 \begin{propos}\label{Pcontrot}
 Let $E$ be an Aubry-Mather set. For every $\varepsilon>0$, there exists a neighborhood $U$ of $E$ that is invariant by the integer translations $(\theta, r)\mapsto (\theta+k, r)$ for $k\in\Z$ and such that every Aubry-Mather set  $\Ec$ that meets $U$ satisfies: $|\rho(E)-\rho(\Ec)|<\varepsilon$.
 \end{propos}

 \demo
 We deduce from Proposition \ref{PAMLip} that $E$ is contained in some strip $\Kc=\R\times [-K, K]$. On such a strip, every $DF^k$ is uniformly bounded.
 
 Let $\Ec$ be an Aubry-Mather set that meets the same strip $\Kc$. Let $(\theta_k, r_k)$ be an orbit in $E$ and $(\alpha_k, \beta_k)$ be an orbit in $\Ec$. We deduce from proposition \ref{PAMhomeo} that for every $k\in \Z$, we have:
 $$|\theta_k-\theta_0+k\rho(E)|\leq 1\quad{\rm and}\quad |\alpha_k-\alpha_0-k\rho(\Ec)|\leq 1.$$
 We deduce
 $$|\rho(E)-\rho(\Ec)|\leq \frac{2}{k}+\frac{|\theta_k-\alpha_k|}{k}+\frac{|\theta_0-\alpha_0|}{k}.$$
 Fixing $k>\frac{4}{\varepsilon}$ large enough, we choose a neighborhood $U$ of $E$ that is invariant by the integer translations $(\theta, r)\mapsto (\theta+k, r)$ for $k\in\Z$ and such that  for every $y=(\alpha, \beta)\in U$, there exists $x=(\theta, r)\in E$ that satisfies $|\theta-\alpha|<\frac{\varepsilon}{4}$ and $\| F^k(x)-F^k(y)\| <\frac{\varepsilon}{4}$. Then for every Aubry-Mather set $\Ec$ that meets $U$, we have  $|\rho(E)-\rho(\Ec)|<\varepsilon$.
 \enddemo
 
  \begin{propos}
Let $F$ be a lift of a conservative twist map $f$. Then for every $\alpha\in\R$, there exists at least one minimizing Aubry-Mather set with rotation number $\alpha$.
 \end{propos}
 
 \demo
 If $\alpha=\frac{p}{q}\in\Q$ is rational, we have proved in Proposition \ref{Pminper} the existence of a $(q,p)$-minimizer $(\theta_k)$. Then the corresponding $F$-orbit $(\theta_k, r_k)$ is minimizing and we deduce from Corollary \ref{Ccomparper} that $E=\{ (\theta_k, r_k)\}+\Z\times\{ 0\}$ is a minimizing Aubry-Mather set with rotation number $\frac{p}{q}$. 
 
 If $\alpha\in\R\backslash\Q$ is irrational, we consider a sequence $(\frac{p_n}{q_n})$ of rational numbers that converge to $\alpha$ and for every $n$ a $(q_n, p_n)$-minimizing orbit $(\theta_k^n, r_k^n)_{k\in\Z}$. As $\theta^n_{q_n}=\theta^n_0+p_n$, there exists $k_n\in [0, p_n-1]$ such that $\theta^n_{k_n+1}-\theta^n_{k_n}\in [0, \frac{p_n}{q_n}]$. Replacing $(\theta_k^n, r_k^n)_{k\in\Z}$ by $(\alpha^n_k, \beta^n_k)=(\theta_{k+k_n}^n-[\theta^n_{k_n}], r_k^n)_{k\in\Z}$ that is also a $(q_n, p_n)$-minimizer, we obtain a sequence of minimizers so that:
 \begin{enumerate}
 \item[$\bullet$] $\alpha_0^n\in [0, 1]$;
 \item[$\bullet$] $(\alpha_1^n-\alpha_0^n)_{n\in\N}$ is bounded and then  $(\alpha_0^n, \beta_0^n)_{n\in\N}$ is also bounded;
 \item[$\bullet$]Êthe rotation number of the $(q_n, p_n)$-minimizer $(\alpha_k^n, \beta_k^n)_{k\in\Z}$ is $\frac{p_n}{q_n}$.
 \end{enumerate}
 We then extract a subsequence so that $(\alpha_0^n, \beta_0^n)_{n\in\N}$ converges to some $(\theta, r)$. Then the orbit of $(\theta, r)$ is also minimizing. If $E={\rm Closure}\left(\{ F^k(\theta, r)+(j, 0); k, j\in\Z\}\right)$, then we deduce from Proposition \ref{PForderedlimit} that $E$ is $F$ ordered and then $E$ is a minimizing Aubry-Mather set.  We deduce from Proposition \ref{Pcontrot} that $\rho(E)=\alpha$.\enddemo

\subsection{Further results on Aubry-Mather sets}
 In \cite{Gol1}, it is proved that the closure of the union of the $\Z\times\{ 0\}$-translated sets of every minimizing orbit is an Aubry-Mather set (hence every minimizing orbit has a rotation number). 
 
In  \cite{Ban}, more precise results concerning the minimizing Aubry-Mather sets are proved. Let us explain them. 
 
We denote the  set of points $(\theta, r)\in\R^2$ having a minimizing orbit by $\Mc(F)$. Then it is closed and   $p(\Mc(F))\subset  \A$ is closed too. The rotation number $\rho:\Mc(F)\rightarrow \R$ is continuous and for every $\alpha\in\R$, the set $\Mc_\alpha(F)=\{ x\in\Mc(F), \rho(x)=\alpha\}$ is non-empty.\\

If $\alpha$ is irrational, then $K_\alpha=p(\Mc_\alpha(F))\subset \A$ is the graph of a Lipschitz map above a compact subset of $\T$. Moreover, there exists a bi-Lipschitz orientation preserving homeomorphims $h:\T\rightarrow \T$ such that 
$$\forall x\in K_\alpha, h(\pi(x))=\pi(f(x)).$$
Hence $K_\alpha$ is:
\begin{enumerate}
\item[-]  either not contained in an invariant loop and then is the union of a Cantor set $C_\alpha$ on which the dynamics is minimal and some homoclinic orbits to $C_\alpha$;
\item[-] or is an invariant loop. In this case the dynamics restricted to $K_\alpha$ can be minimal or Denjoy.
\end{enumerate}
If $\alpha=\frac{p}{q}$ is rational, then $\Mc_\alpha(F)$ is the disjoint union of 3 invariant sets:
\begin{enumerate}
\item[$\bullet$] $\Mc_\alpha^{\rm per}(F)=\{ x\in \Mc_\alpha(F), \pi\circ F^q(x)=\pi(x)+p\}$;
\item[$\bullet$] $\Mc_\alpha^+(F)=\{ x\in \Mc_\alpha(F), \pi\circ F^q(x)>\pi(x)+p\}$;
\item[$\bullet$] $\Mc_\alpha^-(F)=\{ x\in \Mc_\alpha(F), \pi\circ F^q(x)<\pi(x)+p\}$.
\end{enumerate}
The two sets $K_\alpha^+=p(\Mc_\alpha^{\rm per}(F)\cup\Mc_\alpha^+(F))$ and $K_\alpha^-=p(\Mc_\alpha^{\rm per}(F)\cup\Mc_\alpha^-(F))$ are then invariant  Lipschitz graphs above a compact part of $\T$. The points of $p(\Mc_\alpha^+(F)\cup \Mc_\alpha^-(F))$ are   heteroclinic orbits to some periodic points contained in $p(\Mc_\alpha^{\rm per}(F))$.

%Katznelson orntsein?

\section{Ergodic theory for minimizing measures}

\subsection{Green bundles}\label{ssGreen}
We fix a lift $F$ of a conservative twist map. As before $\Mc(F)$ is the set of points whose orbit is minimizing. We use some new notations.

\begin{notas}\label{Nota}{\rm
\begin{enumerate}
\item[$\bullet$] $V(x)=T_x\Vc(x)=\{ 0\}\times\R\subset T_x\R^2$ and for $k\not= 0$, we have $G_k(x)=Df^k(f^{-k}x)V(f^{-k}x)$;
\item[$\bullet$] the slope of $G_k$ (when defined) is denoted by $s_k$: $G_k(x)=\{ (\delta\theta, s_k(x)\delta\theta); \delta\theta\in\R\}$;
\item[$\bullet$] if $\gamma$ is a real Lipschitz function defined on $\T$ or $\R$, then 
$$\gamma'_+(x)=\limsup_{y,z\rightarrow x, y\not=z}\frac{\gamma(y)-\gamma(z)}{y-z}\quad{\rm and}\quad \gamma'_-(t)=\liminf_{y,z\rightarrow x, y\not=z}\frac{\gamma(y)-\gamma(z)}{y-z}.$$
\end{enumerate}}
\end{notas}

We introduce now a set, called ${\rm Green}(f)$. We will see very soon that we can define two natural invariant sub bundles in tangent lines at every point of ${\rm Green}(f)$, that will be very useful in our further study. An important result (see Corollary \ref{CGreenset}) is that all the minimizing  Aubry-Mather sets are contained in ${\rm Green}(f)$.

\begin{notat}{\rm 
We denote by ${\rm Green}(f)$ the set of the points of $\A$ such that {\em along the whole orbit of these points}, we have
$$\forall n\geq 1, s_{-n}(x)<s_{-n-1}(x)<s_{n+1}(x)<s_n(x).$$}
\end{notat} 

\begin{defi}{\rm 
If $x\in {\rm Green}(f)$, the two {\em Green bundles} at $x$ are $G_+(x), G_-(x)\subset T_x(\R^2)$ with slopes $s_-$, $s_+$ where $\displaystyle{s_+(x)=\lim_{n\rightarrow +\infty}s_n(x)}$ and $\displaystyle{ s_-(x)=\lim_{n\rightarrow +\infty}s_{-n}(x)}.$}
\end{defi}
The two Green bundles satisfy the following properties

\begin{propos}\label{PGreen}
Let $f$ be a  conservative twist map. 
\begin{enumerate}
\item[$\bullet$] Then the two Green bundles defined on ${\rm Green}(f)$ are invariant under $Df$: $Df(G_\pm)=G_\pm\circ f$;
\item[$\bullet$] we have $s_+\geq s_-$;
\item[$\bullet$] the map $s_-:\Mc\rightarrow \R$ is lower semi-continuous and the map  $s_+:\Mc\rightarrow \R$ is upper semi-continuous;
\item[$\bullet$] hence $\{ G_-=G_+\}$ is a  $G_\delta$ subset of ${\rm Green}(f)$ and   $s_-=s_+$ is continuous at every point of this set.\end{enumerate} 
\end{propos}

\begin{exer} 
{\rm Prove Proposition \ref{PGreen}.}

\end{exer} 

\begin{thm}\label{TGreen}
Let $f:\A\rightarrow \A$ be a conservative twist map and let $(x_n)_{n\in\Z}$ be the orbit of a point $x=x_0$. The following assertions are equivalent:
\begin{enumerate}
\item[(0)] $x\in{\rm Green}(f)$;
\item[(1)] the projection of every finite segment of the orbit of $x$  is locally minimizing among the segments of points (of $\R$) that have same length and same endpoints;
\item[(2)] along the orbit of $x$, we have for every $k\geq 1$, $s_k>s_{-1}$;
\item[(3)] along the orbit of $x$, we have for every $k\geq 1$, $s_{-k}<s_{1}$;
\item[(4)] there exists a field of half-lines $\delta_+\subset T\A$ along the orbit of $x$ such that:
\begin{enumerate}
\item[$\bullet$] $\delta_+$ is invariant by $Df$: $Df(\delta_+)= \delta_+\circ f$;
\item[$\bullet$] $D\pi\circ \delta_+=\R_+$ ($\delta_+$ is oriented to the right).
\end{enumerate}
\end{enumerate}
\end{thm}

\begin{remas}\label{Rdemid}

{\rm 

\begin{enumerate}
\item Observe that in the point {\em (4)}, you cannot replace `field of half-lines' by `field of lines'. Indeed, along the orbit of every point that is not periodic you can find an invariant field of lines that is transverse to the vertical.
\item in fact, in the proof, we will see that if we denote by $d_+$ the slope of $\delta_+$, we necessarily have $s_-\leq d_+\leq s_+$.
\end{enumerate}

}
\end{remas}

We postpone the proof of Theorem \ref{TGreen} to subsection \ref{ssGreeneq}.
\begin{corol}\label{CGreenset}
Let $f$ be a conservative twist map. Then
\begin{enumerate}
\item[$\bullet$] every accumulation point of an Aubry-Mather set is in ${\rm Green}(f)$;
\item[$\bullet$] every minimizing orbit is in ${\rm Green}(f)$.
\end{enumerate}
\end{corol}
\begin{demo}
$\bullet$ Assume that $x$ is an accumulation point of an invariant Aubry-Mather set $E$. We look at the action of $DF$ on the half-lines $\R_+v$ that are in the tangent space to $\R^2$ along the orbit of $x$. As $E$ is the graph of a Lipschitz map $\gamma:F\rightarrow \R$ and $x$ is an accumulation point of $E$, we have for every $k\in\Z$:
$$\gamma'_+(\pi(F^kx))=\limsup_{y,z\rightarrow \pi(F^kx), y,z\in E,y\not=z}\frac{\gamma(y)-\gamma(z)}{y-z}.$$
This bundle $\Gamma_+=\R_+(1, \gamma'_+)$ in half-lines is transverse to the vertical bundle and invariant by $DF$. We use the characterization {\em (4)} of Theorem \ref{TGreen} to conclude.\\
$\bullet$ The second point of the corollary is a direct consequence of the characterization {\em (1)} of ${\rm Green} (f)$.
\end{demo}

 An interesting consequence of the characterization {\em (1)} of ${\rm Green}(f)$ is
\begin{corol}
The set ${\rm Green}(f)$ is closed.
\end{corol}

When $x$ is a generic point in the support of some hyperbolic measure, $G_-$ is the stable bundle and $G_+$ is the unstable one:

\begin{propos}\label{Pdyncrit} {\bf (Dynamical criterion)} Assume that $x\in{\rm Green}(F)$ has its orbit contained in some strip $\R\times[-K,K]$ (for example $x\in\Mc(F)$ or $x$ is in some Aubry-Mather set) and that $v\in T_x\A$. Then
\begin{enumerate}
\item[$\bullet$] if $\displaystyle{\liminf_{n\rightarrow +\infty} |D(\pi\circ F^n)(x)v|<+\infty}$, then $v\in G_-(x)$;
\item[$\bullet$] if $\displaystyle{\liminf_{n\rightarrow +\infty} |D(\pi\circ F^{-n})(x)v|<+\infty}$, then $v\in G_+(x)$.
\end{enumerate}
\end{propos} 
\demo
We use a symplectic change of linear coordinates along the orbit of $x$ in such a way that the horizontal subspace is now $G_-$ and the vertical subspace doesn't change.

As the orbit of $x$  is contained in some strip $\R\times [-K, K]$,   the slopes $s_{-1}$ and $s_1$ of $G_{-1}=DF^{-1}(V\circ F)$ and $G_{1}=DF(V\circ F^{-1})$ are uniformly bounded along the orbit of $x$. Hence $s_-\in[s_{-1}, s_1]$ is also uniformly bounded and so the changes of basis $P=\begin{pmatrix} 1&0\\ s_-&1
\end{pmatrix}$ as $P^{-1}=\begin{pmatrix} 1&0\\ -s_-&1
\end{pmatrix}$ are also uniformly bounded.
Then the matrix of $DF^n(x)$ in this new basis is
$$\begin{pmatrix}
b_n(x)(s_-(x)-s_{-n}(x))&b_n(x)\\
0&(s_n(F^nx)-s_-(F^nx))b_n(x)
\end{pmatrix}$$
We know that the determinant is $1=(b_n(x))^2(s_-(x)-s_{-n}(x))(s_n(F^nx)-s_-(F^nx))$, that $|s_n(F^nx)-s_-(F^nx)|\leq (s_1(F^nx)-s_{-1}(F^nx))$ is uniformly bounded and that $\displaystyle{ \lim_{n\rightarrow +\infty}(s_-(x)-s_{-n}(x))=0}$; hence $\displaystyle{\lim_{n\rightarrow +\infty}|b_n(x)|=+\infty}$.

Let now $v$ be a vector in $T_x\A$. We denote by $(v_1, v_2)$ its coordinates in the new base we defined just before. Then we have: $ |D(\pi\circ F^n)(x)v|=|b_n(x)|.|(s_-(x)-s_{-n}(x))v_1+v_2|$. As $\displaystyle{ \lim_{n\rightarrow +\infty}(s_-(x)-s_{-n}(x))=0}$ and $\displaystyle{\lim_{n\rightarrow +\infty} |b_n(x)|=+\infty}$, we deduce that if $v_2\not =0$ (i.e. if $v\notin G_-(x)$), then 
 $\displaystyle{\lim_{n\rightarrow +\infty} |D(\pi\circ F^n)(x)v|=+\infty}$.
\enddemo
\begin{exer}\label{Exo}{\rm Let $k\in [1, \infty]$. Let us admit that there exists a dense $G_\delta$ subset $\Gc$ of the set of the $C^k$ conservative twist maps such that for every $f\in\Gc$, for every periodic point $x$ for $f$, if we denote by $N$ the period of $x$, then we have:
\begin{enumerate}
\item[$\bullet$] the eigenvalues of $Df^N(x)$ are distinct;
\item[$\bullet$] every heteroclinic intersection of two hyperbolic periodic orbits is transverse.
\end{enumerate}
Prove that every $f\in \Gc$ has no rational invariant graph.\\
{\bf Hint:} using the  invariance of the Green bundle $G_-$, prove that every periodic point contained in such a rational invariant graph has to be hyperbolic. 
}
\end{exer}

\subsection{Lyapunov exponents and Green bundles}
We have noticed that if a measure $\mu$ with support contained in ${\rm Green}(f)\cap( \R\times[-K, K])$ is hyperbolic, then we have $\mu$-almost everywhere: $E^s=G_-$ and $E^u=G_+$. In this case, we have $G_-\not=G_+$ $\mu$-almost everywhere.\\
We will prove the reverse implication.

\begin{thm}\label{TGreenexp}{\bf (M.-C.~Arnaud)}
Let $f$ be a conservative twist map and let $\mu$ be a measure that is ergodic for $f$, with compact support and such that ${\rm supp}\mu\subset {\rm Green}(f)$. 
Then $d=\dim(G_-\cap G_+)$ is constant $\mu$-almost everywhere and
\begin{enumerate}
\item[$\bullet$] if $d=0$, the measure $\mu$ is hyperbolic with Lyapunov exponents $-\lambda(\mu)<\lambda(\mu)$ given by:
$\lambda(\mu)=\frac{1}{2}\int\log\left( \frac{s_+-s_{-1}}{s_--s_{-1}}\right) d\mu=\frac{1}{2}\int\log\left( 1+\frac{s_+-s_-}{s_--s_{-1}}\right)d\mu$; 
\item[$\bullet$] if $d=1$, the Lyapunov exponents of $\mu$ are zero.
\end{enumerate}
\end{thm}
\begin{rema}{\rm 
Observe that the first part of Theorem \ref{TGreenexp} says to us that   the more distant the Green bundles are, the greater    the positive Lyapunov exponent is.\\
A general result for hyperbolic measures of smooth dynamics is that when the stable and unstable bundles are close together, the Lyapunov exponents are close to zero (see for example \cite{Arna5}).\\
The reverse result is not true in general and what we prove is then specific to the case of the  twist maps. Consider for example the Dirac measure at $(0, 0)$ that is invariant by the linear map of $\R^2$ with matrix $\begin{pmatrix} e^\varepsilon&0\\ 0 & e^{-\varepsilon} \end{pmatrix}$. Then the unstable and stable bundles are $\R\times \{ 0\}$ and $\{ 0\}\times \R$ that are far from each other. But the Lyapunov exponents $\varepsilon$, $-\varepsilon$, can be very close to $0$.}
\end{rema}
\demo
As the dynamics is symplectic, the sum of the Lyapunov exponents is\\
 $\int\log(\det(Df))d\mu=0$, hence there are two Lyapunov exponents $-\lambda(\mu)\leq\lambda(\mu)$. Either these two Lyapunov exponents are zero or the measure is hyperbolic.
\\
We have noticed that when $\mu$ is hyperbolic, then $G_-=E^s\not=G_+=E^u$ $\mu$-almost everywhere. Hence when $d=1$, the Lyapunov exponents are zero. 
Assume now that $d=0$.  Using a bounded change of basis along a generic point for $\mu$ as in the proof of the dynamical criterion, we obtain that $Df(x)_{|G_-(x)}$ is represented by $b_1(x)(s_-(x)-s_{-1}(x))$ and that 
$Df(x)_{|G_+(x)}$ is represented by $b_1(x)(s_+(x)-s_{-1}(x))$. Hence if $v_\pm$ is a base of $G_\pm$, we have:
$$\lambda(v_\pm)=\lim_{n\rightarrow +\infty} \frac{1}{n}\log\left( \|Df^n(x)v_\pm\|\right)=\lim_{n\rightarrow +\infty} \frac{1}{n}\sum_{j=0}^{n-1}\log\left( b_1(f^jx)(s_\pm(f^jx)-s_{-1}(f^jx))\right)
$$
and then by Birkhoff ergodic theorem
$$\lambda(v_+)-\lambda(v_-)= \lim_{n\rightarrow +\infty} \frac{1}{n}\sum_{j=0}^{n-1}\log\left( \frac{s_+(f^jx)-s_{-1}(f^jx)}{s_-(f^jx)-s_{-1}(f^jx)}\right)= \int\log\left( \frac{s_+-s_{-1}}{s_--s_{-1}}\right) d\mu.
$$
As $s_+>s_-$ $\mu$-almost everywhere, we have then $\lambda (v_+)>\lambda(v_-)$. Hence we are in the case of an hyperbolic measure. Then $G_+=E^u$ and $G_-=E^s$ and  $\lambda(v_+)=\lambda(\mu)$, $\lambda(v_-)=-\lambda(\mu)$ and thus $2\lambda(\mu)=  \int\log\left( \frac{s_+-s_{-1}}{s_--s_{-1}}\right) d\mu$.

\enddemo

We have seen in subsection \ref{ssLyapcurves} that the Lyapunov exponents of the measures that are on the irrational invariant curves are zero. But Patrice Le Calvez proved that for general conservative twist maps, many Aubry-Mather sets are (uniformly) hyperbolic, and then are not curves.

\begin{propos} {\bf (P.~Le~Calvez, \cite{LeC2})} Let $k\in [1, \infty]$. There exists a dense $G_\delta$ subset $\Gc_k$ of the set of the $C^k$ conservative twist maps such that for any $f\in \Gc_k$, there exists an open and dense subset $U(f)\subset \R$ such that the minimizing Aubry-Mather sets  having  their rotation number in $U(f)$ are uniformly hyperbolic.

\end{propos}

It may even happen that all the minimizing Aubry-Mather sets are hyperbolic (see \cite{Gor}).

\begin{propos} {\bf (D.~L.~Goroff)} For $|\lambda| >\frac{\sqrt{1+\pi^2}}{\pi}$, the union of the minimizing Aubry-Mather sets for the standard map $f_\lambda$ is uniformly hyperbolic.

\end{propos}
\begin{demo}
We assume that $|\lambda| >\frac{\sqrt{1+\pi^2}}{\pi}$.\\
The standard map with parameter $\lambda$ is defined by $f_\lambda(\theta, r)=(\theta+r+\lambda\sin 2\pi\theta, r+\lambda\sin 2\pi\theta)$ and has the generating function $S_\lambda(\theta, \Theta)=\frac{1}{2}(\Theta-\theta)^2-\frac{\lambda}{2\pi}\cos 2\pi\theta$.\\
Let $E$ be a minimizing  Aubry-Mather set for $f_\lambda$. Observe that $F_\lambda(\theta, r+1)=F_\lambda(\theta, r)+(1, 1)$. Hence we can assume that the rotation number of $E$ is in $(-1, +1)$. Then by the inequalities that we recalled in subsection \ref{sscircle} for circle homeomorphisms, we have for every orbit $(\theta_n, r_n)$ in $E$: $\theta_n-\theta_{n+1}\in (-1, 1)$ and $\theta_n-\theta_{n-1}\in (-1, 1)$ have opposite signs.\\
As $0=\frac{\partial S_\lambda}{\partial \theta}(\theta_n, \theta_{n+1})+\frac{\partial S_\lambda}{\partial\Theta}(\theta_{n-1}, \theta_n)=\theta_n-\theta_{n-1}+\lambda\sin 2\pi\theta_n+\theta_n-\theta_{n-1}$, we deduce that $\lambda\sin 2\pi\theta_n\in (-1, 1)$ i.e. $|\sin 2\pi\theta_n|<\frac{1}{|\lambda|}$. This implies that $|\cos  2\pi\theta_n|>\sqrt{1-\frac{1}{\lambda^2}}$.\\
Moreover, as the orbit is minimizing, we have $$0\leq \frac{\partial^2 S_\lambda}{\partial \theta^2}(\theta_n, \theta_{n+1})+\frac{\partial^2 S_\lambda}{\partial\Theta^2}(\theta_{n-1}, \theta_n)=2+2\pi\lambda\cos 2\pi\theta_n$$ and then $2\geq -2\pi\lambda\cos 2\pi\theta_n$. As $2\pi |\lambda |Ê|\cos 2\pi\theta_n |\geq 2\pi|\lambda | \sqrt{ 1-\frac{1}{\lambda^2}}=2\pi\sqrt{\lambda^2-1}>\frac{2\pi}{\pi}=2$, we have $2\pi\lambda\cos 2\pi\theta_n>0$ and then $2\pi\lambda\cos 2\pi\theta_n>2$.

We can now compute $Df(\theta, r)=\begin{pmatrix} 1+2\pi\lambda \cos 2\pi\theta&1\\ 2\pi\lambda \cos 2\pi\theta & 1\end{pmatrix}$. Observe that $1+2\pi\lambda\cos 2\pi\theta_n>3$ and $2\pi\lambda\cos2\pi\theta_n>2$. Hence if $C=\{ (v_1, v_2)\in\R^2; v_1.v_2\geq 0\}$,   we have $Df(C)\subset C$ and $\forall v\in C, \| Df(v)\|\geq \sqrt{2} \| v\|$ along the orbit $(\theta_n, r_n)$.\\
We have too $(Df(\theta, r))^{-1}=\begin{pmatrix} 1 & -2\pi\lambda\cos 2\pi\theta\\ -1& 1+2\pi\lambda \cos 2\pi\theta\end{pmatrix}$. Hence if $C'=\{ (v_1, v_2)\in\R^2; v_1.v_2\leq 0\}$, we have along the orbit $(\theta_n, r_n)$: $Df^{-1}(C')\subset C'$ and $\forall v\in C', \| Df^{-1}(v)\|\geq \sqrt{2} \| v\|$.

This implies the wanted result.

\end{demo}
A.~Katok proved that the union of the hyperbolic Aubry-Mather sets has zero Lebesgue measure (see \cite{Kat1}). This can be compared to K.A.M. theory that gives in general a union of invariant circles with positive Lebesgue measure.

Theorem \ref{TGreenexp} can be more precise in the case of uniform hyperbolicity.

\begin{propos}\label{PGreenhyp} {\bf (M.-C.~Arnaud)}
Let $M$ be an compact invariant set by a conservative twist map that is contained in ${\rm Green}(f)$. Then $E$ is uniformly hyperbolic if and only if at every point of $M$, the two Green bundles $G_-$ and $G_+$ are transverse.

\end{propos}
We postpone the proof of Proposition \ref{PGreenhyp} to subsection \ref{ssunifhyp}. We don't know if there exist examplesñ  of Aubry-Mather sets that are non-uniformly hyperbolic.

\begin{quest}{\rm 
Does there exist a  conservative twist map that has a non-uniformly hyperbolic Aubry-Mather set?}
\end{quest}

\subsection{Lyapunov exponents and shape on the Aubry-Mather sets}
In the previous subsection, we compared the size of the Lyapunov exponents for the ergodic measures with support in ${\rm Green} (f)$ with the distance between the two Green bundles.  We ask now if we can see a link between the shape of the support of such a measure and the Lyapunov exponents.

\begin{defi}\label{DefC1}{\rm 
Let $M\subset \A$ be a subset of $\A$ and $x\in M$ a   point of $M$. The   {\em paratangent cone} to $M$ at $x$ is the cone of $T_x\A$ denoted by $P_M(x)$ whose elements  are the limits 
$$v=\lim_{ n\rightarrow \infty} \frac{x_n-y_n}{ \tau_n} $$
where $(x_n)$ and $(y_n)$ are sequences of elements of $M$ converging to $x$, $(\tau_n)$ is a sequence of elements of $\R_+^*$ converging to $0$, and $x_n-y_n\in\R^2$, refers to the unique lift of this element of $\A$ that belongs to $[-\frac{1}{2}, \frac{1}{2}[^2$.\\
Here we draw the paratangent cone to a curve at a corner:

 \begin{center}
\includegraphics[width=5cm]{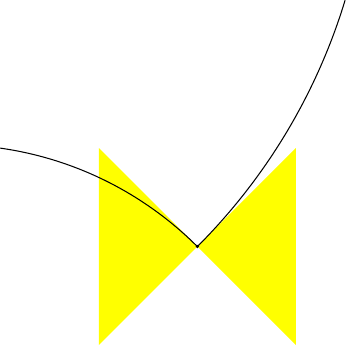}
\end{center}

\noindent We will say that $M$ is {\em $C^1$-regular} at $x$ if there exists a line $D$ of $T_x\A$ such that $P_M(x)\subset D$.\\
If $M$ is not $C^1$-regular at $x$, we say that $M$ is {\em  $C^1$-irregular} at $x$.}
\end{defi}
\begin{rema} {\rm 
Observe that the graph of a Lipschitz map $\gamma$ is $C^1$-regular if and only if $\gamma$ is $C^1$.}
\end{rema}
\begin{notat}{\rm
We denote the set of the slopes of the elements of $P_M(x)$ by $p_M(x)$. }
\end{notat}

\begin{thm}\label{TLyapshape}{\bf (M.-C.~Arnaud)}
Let $\mu$ be an ergodic measure for a conservative twist map with  support in some irrational Aubry-Mather set. Then
\begin{enumerate}
\item[$\bullet$] either the Lyapunov exponents of $\mu$ are zero and ${\rm supp}\mu$ is $C^1$-regular $\mu$-almost everywhere;
\item[$\bullet$] or $\mu$ is hyperbolic and ${\rm supp}\mu$ is $C^1$-irregular $\mu$-almost everywhere.
\end{enumerate}
\end{thm}
\begin{corol}
If the support of an ergodic measure has an irrational rotation number and is contained in some (non necessarily invariant) $C^1$ curve, then its Lyapunov exponents are zero.
\end{corol}
\demo Let $M$ be an irrational Aubry-Mather set and let $\mu$ be the unique ergodic measure with support in $M$. Looking at the proof of Corollary \ref{CGreenset} (see also Remark \ref{Rdemid}), we deduce easily that for $\mu$-almost every point $x\in {\rm supp}\mu$, we have
$$s_-(x)\leq p_M(x)\leq s_+(x).$$
Assume that the Lyapunov exponents of $\mu$ are zero. Then, by theorem \ref{TGreenexp}, we have $\mu$-almost everywhere $G_-=G_+$ i.e. $s_-=s_+$  and then $P_M(x)$ is contained in a line. This exactly means that ${\rm supp}(\mu)$ is $C^1$-regular $\mu$-almost everywhere.  

Now we assume that the Lyapunov exponents of $\mu$ are non zero: $-\lambda(\mu)< \lambda(\mu)$.  The set where ${\rm supp}(\mu)$ is $C^1$-regular is measurable and invariant by $f$. Hence either $\mu$ is $C^1$-regular $\mu$-almost everywhere or $\mu$ is $C^1$-irregular $\mu$-almost everywhere. Assume that $\mu$ is $C^1$-regular $\mu$-almost everywhere.

We will prove the following result (we use for $h'_\pm$ the notations   \ref{Nota}) in subsection  \ref{ssbiLip}.
\begin{propos}\label{PbiLip} {\bf (M.-C.~Arnaud)} Let $h:\T\rightarrow\T$ be a bi-Lipschitz orientation preserving homeomorphism with irrational rotation number. We denote by $\mu$ its unique invariant measure and assume that $h$ is $C^1$-regular $\mu$-almost everywhere. Then uniformly in $\theta\in \T$, we have
$$\lim_{n\rightarrow +\infty} \frac{1}{n}\log\left( h^n\right)'_+=\lim_{n\rightarrow +\infty}\frac{1}{n}\log\left( h^n\right)'_-=0.$$
\end{propos}
Let us explain how we deduce the wanted result. As $M$ is an Aubry-Mather set, it is the graph of a Lipschitz map $\gamma: \pi(M)\rightarrow \R$. We consider the projected-restricted dynamics to $M$, which is $h:\pi(M)\rightarrow\pi(M)$ that is defined by $h(\theta)=\pi\circ f(\theta, \gamma(\theta))$. We denote again by $\mu$ the projected measure $\pi_*\mu$ of $\mu$, that is the unique invariant measure by $h$.  We extend $h$ linearly in its gaps in such a way we obtain a bi-Lipschitz homeomorphism $h$ of $\T$. Because $\mu$ is $C^1$-regular $\mu$-almost everywhere, $h$ is also $C^1$-regular $\mu$-almost everywhere and we deduce from Proposition \ref{PbiLip} that uniformly in $\theta\in \T$, we have
$$\lim_{n\rightarrow +\infty} \frac{1}{n}\log\left( h^n\right)'_+=\lim_{n\rightarrow +\infty}\frac{1}{n}\log\left( h^n\right)'_-=0.$$
Observe that $Df^n(\theta, \gamma(\theta)).(1, \gamma'_+(\theta))=\log\left((h^n)'_+(\theta)\right)(1, \gamma'_+(h^n\theta))$. We deduce that the Lyapunov exponent associated to the vector $(1, \gamma'_+(\theta))$ is zero, which is impossible if the measure is hyperbolic.

\enddemo

 \begin{thm}\label{Tunifhypirr}
  Let $M$ be an irrational Aubry-Mather set of a conservative twist map $f$ of $\A$. Then $M$ is  uniformly hyperbolic  if and only if at every $x\in M$, $M$ is   $C^1$-irregular.
  \end{thm}
  
\demo As $s_-\leq p_M\leq s_+$, if $M$ is $C^1$-irregular everywhere, then $G_-\not= G_+$ at every point of $M$ and by Proposition \ref{PGreenhyp}, $M$ is uniformly hyperbolic.

Assume now that $M$ is uniformly hyperbolic. At first, let us notice that such a $M$ cannot be a curve because of Theorem \ref{Tirrational}.
  
Hence $M$ is a Cantor and the dynamics on $M$ is Lipschitz conjugate to the one of a Denjoy counter-example on its minimal invariant set. Then we consider two points $x\not= y$ of $M$ such that there exists an open  interval $I\subset \T$ whose ends are $\pi(x)$ and $\pi(y)$ and which doesn't meet $\pi(M)$: $I\cap \pi(M)=\emptyset$.  We deduce from the dynamics of the Denjoy counter-examples (see \cite{He1}) that:
\begin{enumerate}
\item[$\bullet$] the positive and negative orbits of $x$ and $y$ under $f$ are dense in $M$;
\item[$\bullet$] $\displaystyle{\lim_{n\rightarrow +\infty}d(f^nx,f^ny)=\lim_{n\rightarrow +\infty}d(f^{-n}x,f^{-n}y)=0}$.
\end{enumerate}
As $M$ is uniformly hyperbolic, we can define a local stable and unstable laminations containing  $M$, $W^s_{\rm loc}$ and $W^u_{\rm loc}$. Then for large  enough $n$, $f^nx$ and $f^ny$ belong  to the same local stable leaf, and $f^{-n}x$ and $f^{-n}y$ belong  to the same local unstable leaf. Hence, because $$\displaystyle{\lim_{n\rightarrow +\infty}d(f^nx,f^ny)=\lim_{n\rightarrow +\infty}d(f^{-n}x,f^{-n}y)=0},$$ for   large enough $n$, the vector joining $f^nx$ to $f^ny$ (resp. $f^{-n}x$ to $f^{-n}y$) is  close  the stable bundle $E^s$ (resp. the unstable bundle $E^u$).

Let now $z\in M$ be any point. Then there exist  two sequences $(i_n)$ and $(j_n)$ of integers which tends to $+\infty$ and are such that:
$$\lim_{n\rightarrow +\infty} f^{i_n}x=\lim_{n\rightarrow +\infty} f^{i_n}y=\lim_{n\rightarrow +\infty} f^{-j_n}x=\lim_{n\rightarrow +\infty} f^{-j_n}y=z.$$
The direction of the ``vector'' joining $f^{i_n}x$ to $f^{i_n}y$ tends to $E^s(z)$ and the direction of the vector joining $f^{-j_n}x$ to $f^{-j_n}y$ tends to $E^u(z)$. Hence: $E^u(z)\cup E^s(z)\subset P_M(z)$ and $M$ is   $C^1$-irregular at $z$.\enddemo

When drawing irrational Aubry-Mather sets that are Cantor sets with the help of a computer, we never observe some angles on these sets. That is why we raise the question:
\begin{quest}{\rm 
Is it possible to draw (with a computer) some irrational Aubry-Mather sets that have some ``corners''?}
\end{quest}
\begin{rema}{\rm 
There is a difficulty in `seing' these corners. On the K.A.M. invariant graphs, the two Green bundles coincide. As $s_+-s_-$ is non-negative and upper-semicontinuous, we deduce that close to the KAM curves, the paratangent cones are very thin, and thus very hard to detect.}
\end{rema}

\section{Complements}
\subsection{Proof of the equivalent definition of a conservative twist map via a generating function}\label{sstwistgene}
We recall the statement.
 \begin{propo}
Let $F:\R^2\rightarrow \R^2$ be a $C^1$ map. Then $F$ is a lift of a conservative twist map $f:¥A\rightarrow \A$  if and only if there exists a $C^2$ function  such that
\begin{enumerate}
\item[$\bullet$] $\forall \theta, \Theta\in\R, S(\theta+1, \Theta +1)=S(\theta, \Theta)$;
\item[$\bullet$] there exists $\varepsilon>0$ so that for all $\theta, \Theta\in\R$, we have 
$$\varepsilon<-\frac{\partial^2S}{\partial \theta \partial \Theta}(\theta, \Theta)<\frac{1}{\varepsilon};$$
\item[$\bullet$] $F(\theta, r)=(\Theta, R)\Longleftrightarrow R=\frac{\partial S}{\partial \Theta}(\theta, \Theta)\quad{\rm and}\quad r=-\frac{\partial S}{\partial \theta}(\theta, \Theta)$.
\end{enumerate}
\end{propo}

\demo $(\Rightarrow)$ Assume that $F:\R^2\rightarrow \R^2$ is the lift of a conservative twist map $f$  such that $\forall x\in \A, \frac{1}{\varepsilon}>D(\pi\circ f)(x)(0, 1)>\varepsilon$.  Then for every $\theta\in\R$, the map $F_\theta:\R\rightarrow \R$ defined by $F_\theta(r)=\pi\circ  F(\theta, r)$ satisfies $\frac{1}{\varepsilon}>F_\theta'>\varepsilon$. Hence every map $F_\theta$ is a $C^1$-diffeomorphism of $\R$ and $G:\R^2\rightarrow \R^2$ defined by $G(\theta, \Theta)=(\theta, F_\theta^{-1}(\Theta))$ is a $C^1$ diffeomorphism.

We introduce the notation $F(\theta, r)=(\Theta(\theta, r), R(\theta, r))$. Note that $G(\theta, \Theta(\theta, r))=(\theta, r)$ i.e. $F_\theta(r)=\Theta(\theta, r)$.
As $f$ is  an exact symplectic twist map, we have: $G^*(f^*\lambda-\lambda)$ is exact. \\
Hence there exists a function $S:\R^2\rightarrow \R$ such that $DS(\theta, \Theta)=R\circ G(\theta, \Theta)d\Theta-F_\theta^{-1}(\Theta)d\theta$. This means exactly that 
$$\frac{\partial S}{\partial \Theta}(\theta, \Theta)=R\circ G(\theta,\Theta)\quad{\rm and}\quad -\frac{\partial S}{\partial \theta}=F_\theta^{-1}(\Theta);$$
and implies that $S$ is $C^2$. Thus we have proved the third point of Proposition \ref{Pgene}.

Let us fix $(\theta, r)\in\A$. We denote by $\gamma$ the loop of $\A$ defined by $\gamma(t)=(\theta+t, r)$ and by $\Gamma$ its lift $\Gamma(t)=(\theta +t, r)$. As $f$ is exact symplectic, we have $\int_\gamma f^*\lambda=\int_\gamma \lambda$. Let us use the notation $F\circ \Gamma(t)=(\Theta_t, R_t)$. As $f$ is isotopic to identity, we have $\Theta_1=\Theta_0+1$.  Moreover: 
$$0=\int_\gamma (f^*\lambda-\lambda)=\int_{G\circ\gamma}G^*(f^*\lambda-\lambda)=\int_{G\circ\gamma} dS=\int_{(\theta+t, \Theta_t)} dS=S(\theta+1, \Theta_0+1)-S(\theta,  \Theta_0);$$
this gives the first point of Proposition \ref{Pgene}.

From $\frac{\partial S}{\partial \theta}(\theta, \Theta(\theta, r))=-r$ we deduce that $\frac{\partial^2 S}{\partial \Theta\partial \theta}(\theta, \Theta(\theta, r)).\frac{\partial \Theta}{\partial r}(\theta, r)=-1$. As $\frac{1}{\varepsilon}>\frac{\partial \Theta}{\partial r}(\theta, r)=D(\pi\circ f)(x)(0, 1)>\varepsilon$, we deduce the second point of Proposition \ref{Pgene}.

($\Leftarrow$) Assume that $S$ satisfies the conclusions of Proposition \ref{Pgene}.

Because of the second point, the maps $\frac{\partial S}{\partial \theta}(\theta, .)$ and $\frac{\partial S}{\partial \Theta}(.,\Theta)$ are $C^1$-diffeomorphisms of $\R$. Hence the third point allows us to define a diffeomorphism $F:\R^2\rightarrow \R^2$. \\
From the first point we deduce that $F(\theta +1, r)=F(\theta, r)+(1, 0)$ hence F is the lift of a $C^1$-diffeomorphism $f:\A\rightarrow \A$.

 Let us prove that $f$ is  a conservative twist map. We use as before the notation $F(\theta, r)=(\Theta(\theta, r), R(\theta, r))$.

From $\frac{\partial S}{\partial \theta}(\theta, \Theta(\theta, r))=-r$ we deduce that $\frac{\partial^2 S}{\partial \Theta\partial \theta}(\theta, \Theta(\theta, r)).\frac{\partial \Theta}{\partial r}(\theta, r)=-1$ and then we have  the twist condition $\varepsilon <\frac{\partial \Theta}{\partial r}(\theta, r)<\frac{1}{\varepsilon}$.

Because $S(\theta+1, \Theta+1)=S(\theta, \Theta)$, we can define a $C^2$-function $s:\A\rightarrow \R$ such that for any lift $\tilde\theta\in\R$ of $\theta$, we have: $s(\theta, r)=S(\tilde\theta, \Theta(\tilde\theta, r))$. Then  $f^*\lambda-\lambda=ds$ is exact. 
In particular, $f$ preserves the orientation. As moreover $F(\theta +1, r)=F(\theta, r)+(1,0)$, we deduce that $f$ is isotopic to identity. Finally, $f$ is conservative.

\enddemo
\subsection{Proof that every invariant continuous graph is minimizing}\label{ssgraphmin}
Let us recall the result due to J.~Mather.
\begin{theo}
 Assume that the graph of a continuous map $\psi:\T\rightarrow \R$ is invariant by a conservative twist map $f$. Then for any generating function associated to $f$, all the orbits contained in the graph of $\psi$ are minimizing.
 \end{theo}
 
 \demo
 Let us introduce the constant $c=\int_0^1\psi(t)dt$ and let us define the $\Z$-periodic $C^1$-function $\eta$ by $\eta(\theta)=\int_0^\theta\psi(t)dt-c\theta$. If $S$ is a generating function of the lift $F$ of $f$ such that $\frac{\partial^2S}{\partial\theta\partial\Theta}<-\varepsilon$, then we define $W(\theta, \Theta)=S(\theta, \Theta)+c(\theta-\Theta)+\eta(\theta)-\eta(\Theta)$.\\
 Observe that $W(\theta+1, \Theta+1)=W(\theta, \Theta)$. Moreover, we have proved in Lemma \ref{Lcoercive} that:
$$\displaystyle{\lim_{|\Theta-\theta|\rightarrow +\infty}\frac{S(\theta, \Theta)}{|\Theta-\theta|}=+\infty
 }.$$
 Hence$\displaystyle{\lim_{|\Theta-\theta|\rightarrow +\infty}\frac{W(\theta, \Theta)}{|\Theta-\theta|}=+\infty
 }$ and hence $W$ has a global minimimum $\mu$. The minimizers of $W$ being critical points, let us look after the critical points of $W$. We have
 $$\frac{\partial W}{\partial \theta} (\theta, \Theta)=\frac{\partial S}{\partial \theta} (\theta, \Theta)+c +\eta'(\theta)= \frac{\partial S}{\partial \theta} (\theta, \Theta)+\psi(\theta);
 $$
 $$\frac{\partial W}{\partial \Theta} (\theta, \Theta)=\frac{\partial S}{\partial \Theta} (\theta, \Theta)-c -\eta'(\Theta)= \frac{\partial S}{\partial \theta} (\theta, \Theta)-\psi(\Theta).
 $$
 Hence $(\theta, \Theta)$ is a critical point if and only $\Theta=\pi\circ F(\theta, \psi(\theta))$. The set of the critical points of $W$ is then a $1$-dimensional connected submanifold of $\R^2$ that corresponds to the graph of $\psi$.  We deduce that the minimum $\mu$ of $W$ is attained exactly on this set. 
 
 Let now $(\theta_k, r_k)_{k\in\Z}$ be the orbit of a point $(\theta, \psi(\theta))$ that is on the invariant graph of $\psi$. Assume that $(\alpha_n)_{\ell\leq n\leq k}$ is a sequence of real numbers so that $\alpha_\ell=\theta_\ell$ and $\alpha_k=\theta_k$. Then  
  $$(k-\ell+1)\mu=\sum_{n=\ell +1}^k W(\theta_{n-1}, \theta_n)  =\sum_{n=\ell +1}^k(S(\theta_{n-1}, \theta_n)+c(\theta_n-\theta_{n-1})+\eta(\theta_{n-1})-\eta(\theta_n))$$
  is less or equal than
  $$\sum_{n=\ell +1}^kW(\alpha_{n-1}, \alpha_n) =\sum_{n=\ell +1}^k(S(\alpha_{n-1}, \alpha_n)+c(\alpha_n-\alpha_{n-1})+\eta(\alpha_{n-1})-\eta(\alpha_n));$$
  i.e. 
  $$\left(\sum_{n=\ell +1}^kS(\theta_{n-1}, \theta_n)\right)+c(\theta_k-\theta_{\ell})+\eta(\theta_{\ell})-\eta(\theta_k)\leq \left(\sum_{n=\ell +1}^kS(\alpha_{n-1}, \alpha_n)\right)+c(\alpha_k-\alpha_{\ell})+\eta(\alpha_{\ell})-\eta(\alpha_k).$$
  As $\alpha_\ell=\theta_\ell$ and $\theta_k=\alpha_k$, we obtain
   $\displaystyle{\sum_{n=\ell +1}^kS(\theta_{n-1}, \theta_n) \leq \sum_{n=\ell +1}^kS(\alpha_{n-1}, \alpha_n)}$   i.e. the orbit of $(\theta, \psi(\theta))$ is minimizing.
 \enddemo
 \subsection{Proof of the equivalence of different definitions of ${\rm Green}(f)$}\label{ssGreeneq}
 The result that we will prove is
 \begin{theo}
 Let $f:\A\rightarrow \A$ be a conservative twist map and let $(x_n)_{n\in\Z}$ be the orbit of a point $x=x_0$. The following assertions are equivalent:
\begin{enumerate}
\item[(0)] $x\in{\rm Green}(f)$;
\item[(1)] the projection of every finite segment of the orbit of $x$  is locally minimizing among the segments of points (of $\R$) that have same length and same endpoints;
\item[(2)] along the orbit of $x$, we have for every $k\geq 1$, $s_k>s_{-1}$;
\item[(3)] along the orbit of $x$, we have for every $k\geq 1$, $s_{-k}<s_{1}$;
\item[(4)] there exists a field of half-lines $\delta_+\subset T\A$ along the orbit of $x$ such that:
\begin{enumerate}
\item[$\bullet$] $\delta_+$ is invariant by $Df$: $Df(\delta_+)= \delta_+\circ f$;
\item[$\bullet$] $D\pi\circ \delta_+=\R_+$ ($\delta_+$ is oriented to the right).
\end{enumerate}
\end{enumerate}
\end{theo}
 
We will use the following notations.

\begin{notas}{\rm
\begin{enumerate}
\item[$\bullet$] $F$ being a lift of $f$, we note:
$$DF^k(y)=\begin{pmatrix}
a_k(y)&b_k(y)\\
c_k(y)&d_k(y)
\end{pmatrix};
$$
\item[$\bullet$] an {\em infinitesimal orbit} along $(x_n)$ is 
$$(\delta\theta_n, \delta r_n)=(Df^n(x)(\delta\theta, \delta r))_{n\in\Z};$$
\item[$\bullet$] a Jacobi field is then the projection  $(\delta\theta_n)_{n\in \N}$ of an infinitesimal orbit;
\item[$\bullet$] if $x_k=(\theta_k, r_k)$, we use the notation 
$$\beta_k=\frac{\partial^2S}{\partial\theta\partial\Theta}(\theta_k, \theta_{k+1}), \quad \alpha_k=\frac{\partial^2S}{\partial \theta^2}(\theta_k, \theta_{k+1})+\frac{\partial^2S}{\partial \Theta^2}(\theta_{k-1}, \theta_{k}).$$
\end{enumerate}}
\end{notas}
\begin{rema}{\rm 
A Jacobi field with  two successive zeroes is the zero field.}
\end{rema}
Let us begin the proof of the theorem.
\medskip

{\bf (1)$\Longrightarrow $(2)} We deduce from the definition of the generating functions that
$$Df(x_k)=\begin{pmatrix}
-\frac{1}{\beta_k}\frac{\partial^2S}{\partial \theta^2}(\theta_k, \theta_{k+1})& -\frac{1}{\beta_k}\\
\beta_k-\frac{1}{\beta_k}\frac{\partial^2S}{\partial\theta^2}(\theta_k, \theta_{k+1})\frac{\partial^2S}{\partial\Theta^2}(\theta_k, \theta_{k+1})& -\frac{1}{\beta_k}\frac{\partial^2S}{\partial\Theta^2}(\theta_k, \theta_{k+1})
\end{pmatrix}.$$ 
Observe too that $(\delta\theta_k)$ is a Jacobi field if and only if for every $k$, we have
$$(*) \beta_{k-1}\delta\theta_{k-1}+\alpha_k\delta\theta_k+\beta_k\delta\theta_{k+1}=0.$$
As we assume that the orbit is locally minimizing, every matrix $H_{n, m}$ is positive semi-definite if:
$$H_{n, m}=\begin{pmatrix}\alpha_{n+1}&\beta_{n+1}&0&.&.&.&0\cr
\beta_{n+1}&\alpha_{n+2}&\beta_{n+2}&.&.&.&0\\
0&\beta_{n+2}&\alpha_{n+3}&.&.&.&0\\
.&.&.&.&.&.&0\\
0&.&.&.&0&\alpha_{m-2}&\beta_{m-2}\\
0&.&.&.&0&\beta_{m-2}&\alpha_{m-1}

\end{pmatrix}
$$
\begin{lemm}
Every matrix $H_{n, m}$ is positive defnite.
\end{lemm}
\begin{demo} Let us assume that $(\delta\theta_k)_{k\in[n+1, m-1]}$ is in the kernel of $H_{n, m}$. Using $(*)$ and the fact that $\beta_k\not=0$ (that is the twist condition), we extend $(\delta\theta_k)$ in a Jacobi field such that $\delta\theta_n=\delta\theta_m=0$.\\
Then, $\delta Q=(0, 0, \delta\theta_{n+1}, \delta\theta_{n+2}, \dots, \delta\theta_{m-2}, \delta\theta_{m-1}, 0, 0)$ is in the isotropic cone of $H_{n-2, m+2}$, and then in its kernel because the matrix is positive semi-definite. Hence  we have a Jacobi field with two successive zeroes, it is the zero field.

\end{demo}
\begin{lemm}
If $k\geq 1$, we have along the orbit of $x$: $s_k>s_{-1}$.
\end{lemm}
\begin{demo}
Let $(\Delta_j)_{j\in[n-k+1, n]}$ be the image by the matrix $H_{n-k, n+1}$ of the Jacobi field $(\delta\theta_j)_{j\in[n-k+1, n]}$ that corresponds to an infinitesimal orbit $(\delta x_j)_{j\in[n-k+1, n]}$ of a vector $\delta x_{n-k}\in V(x_{n-k})$. Then we have
\begin{enumerate}
\item[$\bullet$]  $\Delta_{n-k+1}=0$ because $\delta\theta_{n-k}=0$; 
\item[$\bullet$] for every $j\in[n-k+2, n-2]$, we have $\Delta_j=0$ because we have a Jacobi field;
\item [$\bullet$]  as $\delta x_n=\begin{pmatrix}\delta\theta_n\\ 
s_k(x_n)\delta\theta_n\end{pmatrix}$, we have
$$\Delta_n=\beta_{n-1}\delta\theta_{n-1}+\alpha_n\delta\theta_n=-\beta_n\delta\theta_{n+1}=-\beta_nD(\pi\circ F)\begin{pmatrix} \delta\theta_n\\ s_k(x_n)\delta\theta_n\end{pmatrix}$$
and then
$$\Delta_n=-\beta_n(\beta_n^{-1}(\frac{\partial^2S}{\partial\theta^2}(\theta_n \theta_{n+1})+s_k(x_n)))\delta\theta_n=(s_k(x_n)-s_{-1}(x_n))\delta\theta_n.$$

\end{enumerate}
Finally, we obtain $H_{n-k, n+1}((\delta\theta_j), (\delta\theta_j))=(s_k(x_n)-s_{-1}(x_n))\delta\theta_n^2>0$.
\end{demo}
\medskip

{\bf (2)$\Longrightarrow $(3)}
\begin{lemm}\label{L55}
Assume that we have along the orbit of $x$ and for all $k\geq 1$: $s_k>s_{-1}$. Then we have too along the orbit of $x$: $s_k>s_{k+1}>s_{-1}$.
\end{lemm}
\begin{demo}
We have 
$$Df(x_n)\begin{pmatrix}1\\s_k(x_n)
\end{pmatrix}= \begin{pmatrix}
-\beta_n^{-1}(s_k(x_n)-s_{-1}(x_n))\\ \beta_n-\beta_n^{-1}s_1(x_{n+1})(s_k(x_n)-s_{-1}(x_n))\end{pmatrix}$$
hence $s_{k+1}(x_{n+1})=-\beta_n^2(s_k(x_n)-s_{-1}(x_n))^{-1}+s_1(x_{n+1})$\\
i.e. $(s_{k+1}-s_{-1})(x_{n+1})=(s_1-s_{-1})(x_{n+1})-\beta_n^2(s_k(x_n)-s_{-1}(x_n))^{-1}$\\
and in particular\\
 $(s_{2}-s_{-1})(x_{n+1})=(s_1-s_{-1})(x_{n+1})-\beta_n^2(s_1(x_n)-s_{-1}(x_n))^{-1}$\\
  where $-\beta_n^2(s_1(x_n)-s_{-1}(x_n))^{-1}<0$. Hence $s_2<s_1$.\\
  Substracting what happens for $s_k$ from what happens for $s_{k+1}$ we obtain:
  $$(s_{k+1}-s_k)(x_{n+1})=\beta_n^2\left( (s_{k-1}(x_n)-s_{-1}(x_n))^{-1}-(s_k(x_n)-s_{-1}(x_n))^{-1}\right)$$
  and by recurrence the fact that $(s_k)$ is strictly decreasing.
\end{demo}
\begin{lemm}\label{L56}
If along the orbit of $x$ we have $s_k>s_{k+1}>s_{-1}$ for every $k$, then we have too for every $k$: $s_1>s_{-k}$.
\end{lemm}
\begin{demo}
We assume that $k\geq 2$. We work on the projective space of $\R^2$ that is nothing else than a circle. On this circle, the lines $G_{-1}$, $G_{k+1}$, $G_k$, $G_{k-1}$ are ordered in the direct sense. As $Df^{1-k}$ is symplectic, its projective action preserves the orientation on the circle and then $G_{-k}$, $G_2$, $G_1$ and $V$ are oriented in the direct sense. This means that $s_{-k}<s_2<s_1$.
\end{demo}
\medskip

{\bf (3)$\Longrightarrow $(0)}
Applying results that are analogous to Lemmata \ref{L55} and \ref{L56}, we deduce that if (3) is satisfied, then we have along the orbit of $x$ for every $k\geq 1$: $s_{-1}<s_{k+1}<s_k$ and $s_{-k}<s_{-(k+1)}<s_1$.
\begin{lemm}
Assume that we have $s_1>s_{-k}$ for every $k$ along the orbit of $x$. Then for every $n, k\geq 1$, we have: $s_{-k}<s_n$.
\end{lemm}

\begin{demo}
We assume that $k, n\geq 2$. As in the proof of Lemma \ref{L56}, we work in the projective space. We know that $G_{-1}$, $G_{n+k}$, $G_{n+k-1}$ and $G_{k-1}$ are in the direct sense. Hence their image by $Df^{1-k}$ that are $G_{-k}$, $G_{n+1}$, $G_n$ and $V$ are in the direct sense too, and then $s_{-k}<s_n$.
\end{demo}
\medskip

{\bf (0)$\Longrightarrow $(1)}
We fix a point along the orbit of $x$ (that is denoted by $x$ too) and we go along its orbit until it becomes non strictly
 minimizing. The matrix $H_{0, n}$ is then positive definite but the matrix $H_{0, n+1}$ is not positive definite:
 $$H_{0, n+1}=\begin{pmatrix} \alpha_1&\beta_1&0&\dots&.&0\\
 \beta_1&\alpha_2&.&\dots&.&0\\
 .&\dots&.&\dots&.&0\\
 .&\dots&.&\dots&.&\beta_{n-1}\\
 0&\dots&\dots&0&\beta_{n-1}&\alpha_n\end{pmatrix}.$$
 A vector $(\eta_1, \dots, \eta_n)$ is in the orthogonal subspace to $\R^{n-1}\times\{ 0\}$ for $H_{0,  n+1}$ if and only if we have $\alpha_1\eta_1+\beta_1\eta_2=0$ and for every $j\in[2, n-1]$:
 $\beta_{j-1}\eta_{j-1}+\alpha_j\eta_j+\beta_j\eta_{j+1}=0$, i.e. if $(\eta_j)$ is the projection of an orbit of $V(x)$.
 
Hence if $H_{0, n+1}$ is not positive definite, there exists $\eta_0, \dots, \eta_n$ that is the projection of the orbit of a point of $V(x)\backslash\{ 0\}$ such that:
$$0\geq \eta_n(\beta_{n-1}\eta_{n-1}+\alpha_n\eta_n)=-\beta_n\eta_n\eta_{n+1}.$$
Note that $Df(x_n)=\begin{pmatrix} -b_ns_{-1}&b_n\\ *&*\end{pmatrix}$ hence $\eta_{n+1}=D(\pi\circ f)(x_n)\begin{pmatrix}\eta_n\\ s_n(x_n)\eta_n\end{pmatrix}=b_n(s_n(x_n)-s_{-1}(x_n))\eta_n=-\beta_n^{-1}(s_n(x_n)-s_{-1}(x-n))\eta_n$. We obtain finally $(s_n-s_{-1})(x_n)\eta_n^2\leq 0$. As  $x\in {\rm Green}(f)$, we know that $\eta_n\not=0$. We deduce that $s_n\leq s_{-1}$, a contradiction with the fact that $x\in {\rm Green}(f)$. \\
We deduce that all the matrices $H_{n, m}$ are positive definite and then (1).
\medskip

{\bf (4)$\Longrightarrow $(0)} Now we work on the set of half-lines. We denote by $V_+=\R_+\times\{ 0\}$ the upper vertical and $V_-=-V_+$, $\delta_-=-\delta_+$. This set is a circle and $V_-$, $\delta_+$, $V_+$ and $\delta_-$ are in the direct sense.\\
Because $Df$ preserves the orientation, their images are in the direct sense too, i.e. $\delta_+$, $\R_+(1, s_1)$, $\delta_-$ and $\R_+(-1, -s_1)$ are in the direct sense too. This implies that  $\delta_+$, $\R_+(1, s_1)$, $V_+$, $\delta_-$,  $\R_+(-1, -s_1)$ and $V_-$ are in the direct sense. Taking the images by $Df$, we find that  $\delta_+$, $\R_+(1, s_2)$, $\R_+(1, s_1)$, $\delta_-$,  $\R_+(-1, -s_1)$  and $\R_+(-1, -s_2)$ are in the direct sense and so $\delta_+<s_2<s_1$. Iterating the method, we obtain: $\delta_+<s_{n+1}<s_n$. Replacing $f$ by $f^{-1}$ we obtain too $s_{-n}<s_{-n-1}<\delta_+$.
\medskip

{\bf (0)$\Longrightarrow $(4)} The idea is to use $\delta_+=\R_+(1, s_+)$.

\subsection{Proof of a criterion for uniform hyperbolicity}\label{ssunifhyp}
We want to prove Proposition \ref{PGreenhyp}: 
\begin{propo}  {\bf (M.-C.~Arnaud)}
Let $M$ be a compact invariant set by a conservative twist map that is contained in ${\rm Green}(f)$. Then $M$ is uniformly hyperbolic if and only if at every point of $M$, the two Green bundles $G_-$ and $G_+$ are transverse.
\end{propo}

We have noticed that when $M$ is uniformly hyperbolic, we have $G_-=E^s$ and $G_+=E^u$ on $M$. Hence $G_-$ and $G_+$ are transverse at every point of $M$.

Now we assume that $G_-$ and $G_+$ are transverse at every point of $M$. 

\begin{defi}{\rm
Let $(F_k)_{k\in\Z}$ be a continuous cocycle on a linear normed bundle $P: E\rightarrow K$ above a compact metric space $K$. We say that the cocycle is {\em quasi-hyperbolic} if
$$\forall v\in E, v\not=0\Rightarrow \sup_{k\in\Z}\| F_k v\|=+\infty.$$}
\end{defi}
A consequence of the dynamical criterion (Proposition \ref{Pdyncrit}) is that  if $K\subset \operatorname{Green}  (f)$ is a compact invariant subset of $\operatorname{Green}  (f)$ such that for every $x\in K$, $G_+(x)$ and $G_-(x)$ are transverse, then $(Df^k_{|K})_{k\in\Z}$ is a quasi-hyperbolic cocycle. Hence, we only have to prove the following statement to deduce Proposition \ref{PGreenhyp}.
\begin{thm}\label{Tquasihyp}
Let $(F_k)$ be a continuous, symplectic and quasi-hyperbolic cocycle on a linear and symplectic  (finite dimensional) bundle $P: E\rightarrow K$ above a compact metric space $K$. Then $(F_k)_{k\in\Z}$ is hyperbolic.
\end{thm}

We will deduce Theorem \ref{Tquasihyp} from two lemmata that we will now state and prove. The ideas of the two lemmata and their proofs are   similar  to the ideas  contained in \cite{Man1}  in the setting of the so-called ``quasi-Anosov diffeomorphisms''.

\begin{lemm}\label{Lquasihyp1}
Let $(F_k)_{k\in\Z}$ be a continuous and quasi-hyperbolic cocycle on a linear normed bundle $P: E\rightarrow K$ above a compact metric space $K$. Let us define 
\begin{enumerate}
\item[$\bullet$] $\displaystyle{E^s=\{ v\in E; \sup_{k\geq 0} \| F_kv\|<\infty\}}$;
\item[$\bullet$] $\displaystyle{E^u=\{ v\in E; \sup_{k\leq 0} \| F_kv\|<\infty\}}$.
\end{enumerate}
Then $(F_{n|E^s})_{n\geq 0}$ and $(F_{-n|E^u})_{n\geq 0}$ are uniformly contracting.
\end{lemm} 

\begin{lemm}\label{Lquasihyp2}
Let $(F_k)_{k\in\Z}$ be a continuous and quasi-hyperbolic cocycle on a linear normed bundle $P: E\rightarrow K$ above a compact metric space $K$. We denote by $f_k:K\rightarrow K$ the underlying dynamics such that $f_k\circ P=P\circ F_k$. If $(x_n)$ is a sequence of points of $K$ tending to $x$ and $(k_n)$ a sequence of integers tending to $+\infty$ such that $\displaystyle{\lim_{n\rightarrow \infty}   f_{k_n}(x_n)=y\in K}$, then $\dim E^u(y)\geq {\rm codim} E^s(x)$.
\end{lemm}

\noindent Let us explain how to deduce Theorem \ref{Tquasihyp} from these lemmata:
\smallskip

\noindent{\bf Proof of theorem \ref{Tquasihyp}:} If the dimension of $E$ is $2d$, we only have to prove that: $\forall x\in K, \dim E^u(x)=\dim E^s(x)=d$. Let us prove for example that $\dim E^u(x)=d$. \\
By lemma \ref{Lquasihyp1}, $(F_{n|E^s})_{n\geq 0}$ and $(F_{-n|E^u})_{n\geq 0}$ are uniformly contracting. As the cocycle is symplectic, we deduce that every $E^s(x)$ and $E^u(x)$ is isotropic for the symplectic form  and then $\dim E^s(x)\leq d$ and $\dim E^u(x)\leq d$.\\
 Let us now consider $x\in K$. As $K$ is compact, we can find a sequence $(k_n)_{n\in\N}$ of integers tending to $+\infty$ such that the sequence $(f_{k_n}(x))_{n\in\N}$ converges to a point $y\in K$. Then, by Lemma \ref{Lquasihyp2}, we have:  $\dim E^u(y)\geq {\rm codim} E^s(x)$. But we know that $\dim E^u(y)\leq d$, hence  $  2d-\dim E^s(x)\leq \dim E^u(y)\leq d$ and $\dim E^s(x)=d$.
 \enddemo
 Let us now prove the two lemmata.
 
 \noindent{\bf Proof of lemma \ref{Lquasihyp1}:} We will only prove the result for $E^s$.
 
 Let us assume that we know that:
 $$(*)\quad \forall C>1, \exists N_C\geq 1, \forall v\in E^s, \forall n\geq N_C, \| F_nv\| \leq \frac{ \sup\{ \| F_kv\| ; k\geq 0\}}{C}.$$
 We choose $C>1$. Then  $\sup \{ \| F_k v\| ; k\geq 0\} =\sup\{ \| F_k v\| ; k\in |[0, N_C]|\}$. We define: $M=\sup\{ \| F_k (x)\| ; x\in K, k\in |[0, N_C]|\}$.   Then, if $j\in |[0, N_C-1]|$ and $n\in\N$:
 $$\| F_{nN_c+j}v\| \leq \frac{1}{C} \sup\{ \|  F_{(n-1)N_C+j+k}v\|; k\geq 0\}\leq \frac{1}{C^2} \sup\{ \|  F_{(n-2)N_C+j+k}v\|; k\geq 0\}$$
 $$\dots \leq \frac{1}{C^n} \sup\{ \|  F_{j+k}v\|; k\geq 0\}\leq \frac{1}{C^n} \sup\{ \|  F_{k}v\|; k\geq 0\}\leq \frac{M}{C^n}\| v\|.$$
 This prove exponential contraction. \\
 
 Let us now prove $(*)$. If $(*)$ is not true, there exists $C>1$, a sequence $(k_n)$ in $\N$ tending to $+\infty$ and $v_n\in E^s$ with $\| v_n\|=1$ such that:
 $$\forall n\in\N, \| F_{k_n}v_n\| \geq \frac{\sup\{ \| F_k v_n\|; k\geq 0\}}{C}.$$ 
We define: $w_n=\frac{F_{k_n}(v_n)}{\| F_{k_n}(v_n)\|}$. Taking a subsequence, we can assume that the sequence $(w_n)$ converges to a limit $w\in E$. Then we have:
 $$\forall n\in\N, \forall k\in[-k_n, +\infty[, \| F_kw_n\|=\frac{\| F_{k+k_n}(v_n)\|}{\| F_{k_n}v_n\|}\leq \frac{\sup\{ \| F_jv_n\|; j\geq 0\} }{\| F_{k_n}v_n\|}\leq C.$$
 Hence, $\forall k\in\Z, \| F_k w\|\leq C$. This is impossible because $\| w\|=1$ and the cocycle is quasi-hyperbolic.\enddemo

  \noindent{\bf Proof of lemma \ref{Lquasihyp2}:} With the notation of this lemma, we choose a linear subspace $V\subset E_x$ such that $V$ is transverse to $E^s(x)$. We want to prove that $\dim E^u(y)\geq \dim V$.\\
  We choose $V_n\subset  E_{x_n}$ such that $\displaystyle{\lim_{n\rightarrow \infty} V_n=V}$. Extracting a subsequence, we have: $\displaystyle{\lim_{n\rightarrow \infty} F_{k_n}(V_n)=V'\subset E_y}$. Then  we will prove that $V'\subset E^u(y)$.\\
  
  Let us assume that we have proved that there exists $C>0$ such that 
  $$(*)\quad  \forall n, \forall 0\leq k\leq k_n,  \| F_{-k|F_{k_n}(V_n)}\|\leq C.$$
  Then,  $\forall w\in V', \forall k\in \Z_-, \| F_kw\| \leq C\| w\|$ and $w\in E^u(y)$.
  
  Let us now assume that $(*)$ is not true. Replacing $(k_n)$ by a subsequence, we find for all $n\in N$ an integer $i_n$ between $0$ and $k_n$ such that     $\| F_{-i_n|F_{k_n}(V_n)}\|\geq n$. We choose $w_n\in F_{k_n}(V_n)$ such that $\| w_n\| =1$ and $\|F_{-i_n}(w_n)\| =\|  F_{-i_n|F_{k_n}(V_n)}\|$. We may even assume that: $\| F_{-i_n}(w_n)\| = \sup\{ \| F_k(w_n)\| ; -k_n\leq k\leq 0\} \geq n$.\\
  Then $\displaystyle{\lim_{n\rightarrow +\infty} i_n=+\infty}$. If $v_n=\frac{F_{-i_n}(w_n)}{\|F_{-i_n}(w_n)\|}$, we may extract a subsequence and assume that: $\displaystyle{\lim_{n\rightarrow \infty} v_n=v}$, with $\| v\|=1$.\\
   Then we have 
  $\forall k\in |[ 0, i_n]|, \| F_kv_n\| \leq \| v_n\|$ for 	all $ k= 0, \dots,  i_n$, and therefore  $ \| F_k v\|\leq \| v\|$ for all $k\in\N$ and $v\in E^s$.\\
Now, we have two cases:
  \begin{enumerate}
  \item[$\bullet$] either $(k_n-i_n)$ doesn't tend to $+\infty$; we may extract a subsequence and assume that $\displaystyle{\lim_{n\rightarrow +\infty}(k_n-i_n)=N\geq 0}$; then: $\displaystyle{  F_{-N}v=\lim_{n\rightarrow \infty}F_{i_n-k_n}(v_n)=\lim_{n\rightarrow \infty}\frac{F_{-k_n}(w_n)}{\|F_{-i_n}(w_n)\|}}$. We have: $\frac{F_{-k_n}(w_n)}{\|F_{-i_n}(w_n)\|}\in V_n$ and then $F_{-N}v\in V$. Moreover, $F_{-N}v\in F_{-N}E^s=E^s$. 
  As $\| v\| =1$ and $V$ is transverse to $E^s_x$, we obtain a contradiction.
  \item[$\bullet$] or $\displaystyle{\lim_{n\rightarrow \infty} (k_n-i_n)=+\infty}$. In this case, for every $k=-k_n+i_n, \dots, i_n$, we have $-k_n\leq k-i_n\leq 0$ and therefore $\| F_kv_n\|=\frac{\| F_{k-i_n}w_n\|}{\| F_{-i_n}w_n\|}\leq 1=\| v_n\|$. Hence, since $v_n\rightarrow v$, $i_n\rightarrow +\infty$, and $-k_n+i_n\rightarrow -\infty$, when $n\rightarrow +\infty$, we obtain $\| F_kv\|\leq \| v\|=1$, for all $k\in\Z$. This implies
  $v\in E^s\cap E^u$. This contradicts $\| v\|=1$ and the fact that the cocycle is quasi-hyperbolic.
  \end{enumerate}
  \enddemo
  \subsection{Proof of Proposition \ref{PbiLip}} \label{ssbiLip}
  We will prove
  \begin{propo} {\bf (M.-C.~Arnaud)} Let $h:\T\rightarrow\T$ be a bi-Lipschitz orientation preserving homeomorphism with irrational rotation number. We denote by $\mu$ its unique invariant measure and assume that $h$ is $C^1$-regular $\mu$-almost everywhere. Then uniformly in $\theta\in \T$, we have
$$\lim_{n\rightarrow +\infty} \frac{1}{n}\log\left( h^n\right)'_+=\lim_{n\rightarrow +\infty}\frac{1}{n}\log\left( h^n\right)'_-=0.$$
\end{propo}

\demo
A fundamental argument of the proof is a result proved by A.~Furman in \cite{Fur} that is an improvement of Kingman subadditive theorem in the case of a  unique ergodic measure.

\begin{thm}\label{Tfur}{\bf (A.~Furman)}
Let $(X,\mu)$ be a Borel probability space, $T$ be a continuous measure preserving transformation of $(X,\mu)$  such that $\mu$ is uniquely ergodic for $T$    and llet $(f_n)\in L^1(X,\mu )$ be a $T$-sub-additive sequence of upper semi-continuous functions.  Let $\displaystyle{\Lambda ((f_k))=\lim_{n\rightarrow \infty}\frac{1}{n}\int f_nd\mu}$ be the constant associated to $f$ via the sub-additive ergodic theorem. Then:
$$\forall \varepsilon >0, \exists N\geq 0, \forall n\geq N, \forall x\in X, \frac{1}{n}f_n(x)\leq \Lambda ((f_k))+\varepsilon.$$
\end{thm}
We apply Theorem \ref{Tfur} for $(X, \mu)= (\T, \mu)$, $T=h$   (resp. $T=h^{-1}$) and $f_n=-\log\left( (h^n)'_-\right)$  (resp. $f_n=-\log\left((h^{-n})'_-\right)$). Fixing $\varepsilon>0$, we find $N\geq 0$ such that for every $n\geq N$ and every $\theta\in\T$, we have
$$-\frac{1}{n}\log\left( (h^n)'_-(\theta)\right)\leq \Lambda((f_k))+\varepsilon.$$
We denote by $d\theta$ the Lebesgue measure on $\T$. Because of Jensen inequality for the convex function $-\log$, we have $$-\log\left( \int\left( (h^n)'_-\right)d\theta\right) \leq-\int\log\left( (h^n)'_-\right)d\theta.$$
Moreover, if $H$ is a lift of $h$, 
$$\int\left( (h^n)'_-\right)d\theta\leq \int\left(  h^n \right)'d\theta=\left[ H^n\right]_0^1=1.$$
We deduce
$$\Lambda ((f_k))+\varepsilon\geq -\frac{1}{n}\int\log\left(  h^n \right)'_-d\theta\geq -\log 1=0$$
and then $\Lambda((f_k))\geq 0$.

Finally, we obtain in particular: $\Lambda \left(-\log ( (h^n)'_-)\right)\geq 0$ and $\Lambda \left(-\log ( (h^{-n})'_-)\right)\geq 0$. \\
observe that $\left(h^{-n}\right)'_-(\theta)=\frac{1}{\left( h^n\right)'_+(h^{-n}\theta)}$ hence 
$$\int\log\left(h^{-n}\right)'_- d\mu= -\int\log\left( h^n\right)'_+d\mu.$$
Because $h$ is $C^1$-regular $\mu$-almost everywhere we have $\mu$-almost everywhere $\displaystyle{\prod_{j=0}^{n-1}h'_+(h^j\theta)=\prod_{j=0}^{n-1}h'_-(h^j\theta)}$. Because $(h^n)'_-$ and $(h^n)'_+$ are between these two numbers, we deduce that we have $\mu$-almost everywhere   $(h^n)'_-(\theta)=  (h^n)'_+(\theta)$ and then
 $$\frac{1}{n}\int\log\left(h^{-n}\right)'_- d\mu= -\frac{1}{n}\int\log\left( h^n\right)'_-d\mu$$
 and
 $$\Lambda \left(-\log ( (h^n)'_-)\right)=-\Lambda \left(-\log ( (h^{-n})'_-)\right)= 0.$$
 We deduce then from Theorem \ref{Tfur} that for every $\varepsilon>0$, there exists $N\geq 0$ such that for every $n\geq N$ and every $\theta\in\T$, we have
 $$-\frac{1}{n}\log\left( (h^n)'_-(\theta)\right)\leq \varepsilon\quad{\rm and}\quad  \frac{1}{n}\log\left( (h^n)'_+\theta)\right)=-\frac{1}{n}\log\left( (h^n)'_-(h^n\theta)\right)\leq  \varepsilon$$
 then 
 $$-\varepsilon \leq \frac{1}{n}\log\left( (h^n)'_-(\theta)\right)\leq \frac{1}{n}\log\left( (h^n)'_+\theta)\right)\leq\varepsilon.$$

\enddemo

\end{document}